\PassOptionsToPackage{left=25mm,right=25mm,top=27mm,bottom=27mm}{geometry}
\documentclass[preprint,11pt,authoryear]{elsarticle}

\usepackage{geometry}
\usepackage[T1]{fontenc}
\usepackage[utf8]{inputenc}
\usepackage{amsmath,amssymb,amsfonts,bm}
\usepackage{booktabs,array,tabularx,multirow}
\usepackage{longtable}
\newcolumntype{L}[1]{>{\raggedright\arraybackslash}p{#1}}
\newcolumntype{Y}{>{\raggedright\arraybackslash}X}
\usepackage{graphicx}
\usepackage{eso-pic}
\usepackage{float}
\usepackage{flafter}
\usepackage[section]{placeins}
\usepackage[skip=5pt]{caption}
\usepackage{subcaption}
\usepackage{tikz}
\usetikzlibrary{arrows.meta,positioning,fit,calc}
\usepackage{xcolor}
\definecolor{darkred}{RGB}{139,0,0}
\definecolor{linkblue}{RGB}{0,82,155}
\usepackage{xurl}
\usepackage[
  colorlinks=true,
  linkcolor=linkblue,
  citecolor=linkblue,
  urlcolor=linkblue,
  breaklinks=true
]{hyperref}

\usepackage{natbib}
\usepackage{enumitem}
\usepackage{setspace}
\setlist[itemize]{leftmargin=*,nosep}
\graphicspath{{figures/}}

\begin{document}
\AddToShipoutPictureFG*{%
  \AtPageUpperLeft{%
    \raisebox{-1.60cm}[0pt][0pt]{%
      \hspace*{\dimexpr 1in+\oddsidemargin\relax}%
      \includegraphics[height=1.3cm]{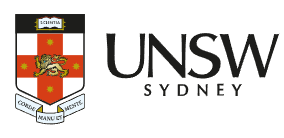}%
      \hspace{0.65cm}%
      \includegraphics[height=1.2cm]{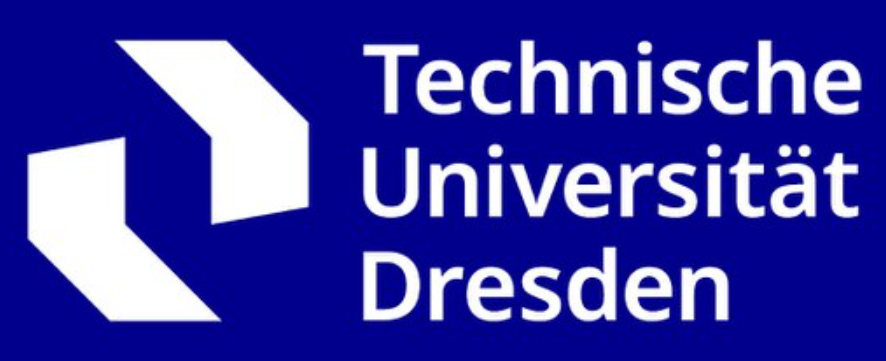}%
      \hspace{0.65cm}%
      \includegraphics[height=1.3cm]{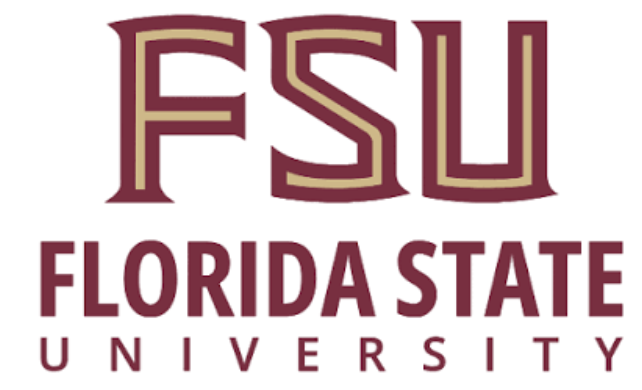}%
      \hspace{0.65cm}%
      \includegraphics[height=1.3cm]{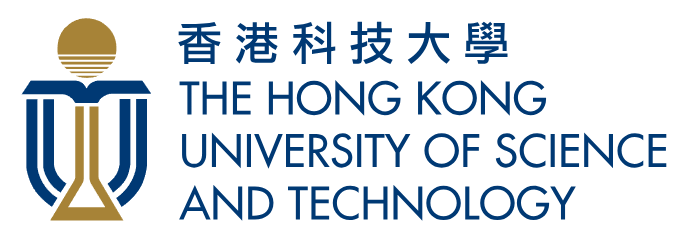}%
      \hspace{0.65cm}%
      \includegraphics[height=1.3cm]{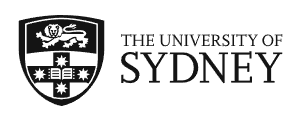}%
      \hspace{0.65cm}%
    }%
  }%
  \AtPageUpperLeft{%
    \raisebox{-2.05cm}[0pt][0pt]{%
      \hspace*{\dimexpr 1in+\oddsidemargin\relax}%
      \rule{\textwidth}{0.4pt}%
    }%
  }%
}

\begin{frontmatter}

\title{{\bfseries
{\LARGE ResiliFlow:} An Open Transport World Model for Infrastructure Perception
and Disaster Resilience}}

\author[aff1]{Junxiang Xu\corref{cor1}}
\ead{junxiang.xu@unsw.edu.au}

\author[aff1]{Vinayak Dixit}

\author[aff1,aff2]{S. Travis Waller}

\author[aff1]{Divya Jayakumar Nair}

\author[aff3]{Qianwen (Vivian) Guo}

\author[aff4]{Sisi Jian}

\author[aff4]{Xiao Wen}

\author[aff1]{Ashutosh Ashutosh}

\author[aff1]{Sunhyung Yoo}

\author[aff1]{Julius Secadiningrat}

\author[aff5]{Jingni Guo}

\cortext[cor1]{Corresponding author.}

\affiliation[aff1]{
  organization={Research Centre for Integrated Transport Innovation (rCITI), School of Civil and Environmental Engineering, UNSW Sydney},
  addressline={Kensington},
  city={Sydney},
  state={NSW},
  postcode={2052},
  country={Australia}
}

\affiliation[aff2]{
  organization={Lighthouse Professorship ``Transport Modelling and Simulation'', Faculty of Transport and Traffic Sciences, Technische Universität Dresden},
  city={Dresden},
  country={Germany}
}

\affiliation[aff3]{
  organization={Department of Civil and Environmental Engineering, Florida State University},
  city={Tallahassee},
  state={FL},
  postcode={32310},
  country={USA}
}

\affiliation[aff4]{
  organization={Department of Civil and Environmental Engineering, The Hong Kong University of Science and Technology},
  city={Hong Kong},
  country={China}
}

\affiliation[aff5]{
  organization={Institute of Transport and Logistics Studies, School of International Business, Logistics and Supply Chains, The University of Sydney},
  city={Sydney},
  state={NSW},
  postcode={2006},
  country={Australia}
}

\begin{abstract}
Transport resilience work is often split across separate data preparation scripts, network models, simulation tools, image inspection systems and reports. This fragmentation makes it difficult to move from an observation to a tested and reviewable decision. We introduce ResiliFlow, an open transport world model concept and an implemented platform for infrastructure resilience, response and recovery. The platform connects two workspaces. Disaster Transport Resilience Analysis provides six map-centred functions for critical-road and critical-area identification, recovery prioritisation, disruption routing, resilience testing and scenario simulation. AI-based Transport Infrastructure Perception and Decision Support organises street-level and satellite evidence, detects visible road, footpath and kerb conditions, and prepares these observations for human-reviewed intervention planning. Both workspaces share an eight-step cycle of perception, prediction, model development, verification, execution, decision, feedback and memory. Research Validation records assumptions and checks, while a local Assistant and an optional multi-provider large language model Copilot translate user questions into bounded calls to executable tools. We document the platform architecture, representative mathematical models, interface evidence and computer-vision learning results. Examples show accurate recognition across eight visible-condition classes, while compact error analysis demonstrates how difficult cases guide continued learning. ResiliFlow shows how transport models can become an inspectable, reusable and question-led system rather than a collection of disconnected analyses. The accompanying release is intended to support research collaboration, public scrutiny and extension under institutional review.
\end{abstract}


\end{frontmatter}

\section{Introduction}

Transport agencies rarely lack models. They often lack a practical way to connect evidence, models and decisions. A flood depth layer may sit in one system, a road network in another, an evacuation model in a third and an inspection image archive in a fourth. Analysts then move files between scripts, maps and reports. Each transfer can lose assumptions, provenance and the reason why a result was produced. The problem becomes more visible during a disruption, when an operator needs to know what changed, why it matters, what can be done and how confident the answer should be.

Transport resilience research provides strong foundations for answering parts of this chain. It has established ways to describe system performance, vulnerability, recovery and community resilience \citep{bruneau2003a,cutter2008a,mattsson2015vulnerability,gonalves2020resilience,bergantino2024assessing}. Flood impact functions link physical exposure to road disruption \citep{pregnolato2017the}. Network science explains how topology and efficiency shape failure consequences \citep{latora2001efficient}. Evacuation research connects route choice, traffic dynamics, shelter assignment and operational response \citep{dixit2009hurricane,wolshon2010stateoftheart,ng2010reliable,ng2010a,dixit2011validation,dixit2012modeling,montz2013assessing,dixit2014evacuation}. These contributions are valuable, yet their outputs are usually delivered through separate studies and specialised interfaces.

Recent work has widened this evidence base. \citet{niu2022linklevel} used crowd-sourced data for link-level resilience analysis. \citet{niu2022are} examined wildfire fatalities in relation to road-network characteristics. \citet{yazdani2021hospital} reviewed hospital evacuation modelling, and \citet{yazdani2022an} developed an integrated flood-evacuation decision model. \citet{niu2023predisaster} addressed resilient road investment under uncertainty. \citet{xu2024exploring} compared five traffic equilibrium formulations, while \citet{xu2024predisaster} and \citet{xu2025predisaster} studied evacuation network design under uncertain demand, behavioural choice and connectivity reliability. \citet{xu2026if} added stated-preference evidence on flood evacuation behaviour. This progression shows that transport resilience decisions need behavioural, network, spatial and operational evidence in the same workflow.

Artificial intelligence extends what can be observed. Road-damage datasets and modern detectors can locate visible defects in street-level imagery \citep{arya2024rdd,li2024rddyolo,zhao2024detrs,jocher2026yolo26}. Urban-scene datasets support broader perception research \citep{cordts2016the,yu2020bddk}. Machine learning has also been applied to bridge cracks, disaster imagery, remote sensing and emergency planning \citep{munawar2022modern,munawar2022disaster,munawar2022remote,munawar2022an}. However, a defect detector usually ends at a bounding box. It does not automatically explain whether a footpath repair reconnects a school to public transport, whether a road closure changes regional accessibility, or whether the evidence is suitable for an operational claim.

This paper introduces ResiliFlow as a step towards that connected capability. We define a \emph{transport world model} as an executable and inspectable system that maintains a shared transport state, selects and runs appropriate models, records evidence, compares actions and updates its state through feedback. The term does not imply a perfect replica of the world. It describes a disciplined cycle that links perception, prediction, model development, verification, execution, decision, feedback and memory. The first implemented platform brings together two workspaces:

\begin{itemize}
  \item \textbf{Disaster Transport Resilience Analysis}, which provides six connected, map-centred functions for disruption diagnosis, routing, resilience testing, recovery planning and scenario simulation.
  \item \textbf{AI-based Transport Infrastructure Perception and Decision Support}, which organises lawful imagery, detects visible road, footpath and kerb conditions, records uncertainty and prepares observations for intervention planning.
\end{itemize}

Table~\ref{tab:workspaces} shows how these workspaces serve different starting points while contributing to one transport state.

\begin{table}[H]
\centering
\caption{The two ResiliFlow workspaces and their immediate user value.}
\label{tab:workspaces}
\footnotesize
\begin{tabularx}{\textwidth}{L{3.6cm}L{3.2cm}YY}
\toprule
Workspace & Typical starting point & What the user can do & Evidence returned \\
\midrule
Disaster Transport Resilience Analysis
& Network, disruption or scenario data
& Locate critical roads and areas, test alternatives, plan recovery and run scenarios
& Maps, rankings, routes, resilience curves, simulations and validation records \\
AI-based Transport Infrastructure Perception and Decision Support
& A place, road view or authorised image
& Inspect visible road, footpath and kerb conditions and move from observations towards repair priorities
& Localised conditions, confidence, provenance, review state and decision context \\
\bottomrule
\end{tabularx}
\end{table}
\FloatBarrier

The contribution is fourfold. First, the paper provides a readable operational definition of a transport world model. Second, it presents an implemented architecture in which data, model outputs, validation traces and memory remain linked. Third, it shows how local guidance and optional large language model support can translate questions into bounded calls to executable tools. Fourth, it demonstrates how learned infrastructure perception can lead to mapped conditions and decision support. The aim is to invite collaboration around an ambitious, transparent platform idea without publishing security-sensitive implementation details.

\section{Related foundations and the remaining gap}

\subsection{Disaster transport resilience and recovery}

Resilience is broader than robustness. Robustness asks how well a network tolerates disruption, while resilience also concerns response, restoration and adaptation. Quantitative frameworks have used functionality curves, redundancy, accessibility, efficiency and recovery time \citep{bruneau2003a,mattsson2015vulnerability,gonalves2020resilience}. Empirical reviews have increasingly emphasised observations from real systems \citep{bergantino2024assessing}. At the same time, place-based research has shown that infrastructure performance and community consequences cannot be separated \citep{cutter2008a,rottemberg2022inequality}.

The research underpinning ResiliFlow spans several parts of this problem. \citet{xu2024datadriven} demonstrated the importance of treating connectivity as a measured network property, and \citet{xu2025reevaluating} examined the role of the relative size of the largest connected component. \citet{waller2025rapid} developed a rapid post-disruption assessment under limited information. \citet{zhang2025integrating} linked geographical equity and travel behaviour to network resilience enhancement. \citet{zhang2023equity} studied equity at the pre-event optimisation stage, while \citet{najmi2023equity} formalised equity in network design and pricing. \citet{xu2025implementing} discussed equitable wildfire response plans. Together, these studies motivate a workspace that keeps topology, behaviour, uncertainty, equity and recovery priorities visible within the same project.

Risk also changes with interaction between networks and decisions. \citet{guo2021risk} assessed multimodal transport-network risk, and \citet{guo2020risk} studied uncertain cascading failure. \citet{xu2026ethical} and \citet{xu2026from} showed why disaster decisions can create ethical cascades when one apparently reasonable choice transfers risk to another group or location. \citet{xu2025artificial} reviewed achievements and challenges for artificial intelligence in disaster management. ResiliFlow responds by retaining model assumptions, evidence sources, uncertainty and human review status alongside each result.

\subsection{Digital twins and world models}

Digital twins connect physical systems with digital representations and are increasingly used in cities and transport \citep{batty2018digital,dembski2020urban,wang2024architecture}. Transport-oriented studies have described digital-twin architectures for adaptability and resilience \citep{feng2023resilience,wang2024digital}. Disaster researchers have proposed urban twins for risk management and post-disaster coordination \citep{macatulad2024continuing,lagap2024digital}. These systems provide important data integration and monitoring capabilities.

ResiliFlow places a different emphasis on executable reasoning. Its core object is not only a synchronised representation. It is a traceable state transition that connects a user question to model selection, verification, execution, comparison and reuse. Recent machine-learning world models have shown how learned internal states can support diverse control tasks \citep{hafner2025mastering}. The transport setting needs additional safeguards because many decisions depend on established network models, statutory data, human judgement and reproducible calculations. ResiliFlow therefore combines learned perception with explicit models and validation gates. This hybrid design makes the reason for an output reviewable.

\subsection{AI inspection and trustworthy evidence}

Street-level perception has advanced through larger datasets and real-time object detectors. RDD2022 assembled multinational road-damage imagery \citep{arya2024rdd}. RDD-YOLO adapted a modern one-stage detector to road defects \citep{li2024rddyolo}. RT-DETR demonstrated an end-to-end transformer detector with real-time performance \citep{zhao2024detrs}. YOLO26 extended real-time detection with an end-to-end family that includes small and larger variants \citep{jocher2026yolo26}. Cityscapes and BDD100K established broader urban-scene benchmarks \citep{cordts2016the,yu2020bddk}, while unpaired image translation offered one path for appearance adaptation across domains \citep{zhu2017unpaired}.

High benchmark scores do not establish safe operation in a new city. Camera position, road material, markings, drainage, maintenance practice and weather can all shift. A detector may recognise an easy close-up crack but miss several small defects in a wide image. It may also mistake a manhole for a pothole. For this reason, the perception workspace uses dataset records, model cards, provenance and independent test gates \citep{mitchell2019model,gebru2021datasheets,wilkinson2016the}. Reporting guidance from other high-stakes machine-learning fields also reinforces the value of explicit data, evaluation and intended-use statements \citep{collins2024tripodai}. The design principle is simple: confidence is useful only when its meaning has been tested.

\subsection{LLM interfaces and reproducible workflows}

Large language models can make complex modelling systems easier to navigate. Retrieval-augmented generation connects language output to an external evidence store \citep{lewis2020rag}. ReAct connects reasoning with actions \citep{yao2023react}. Toolformer showed that language models can learn when and how to call external tools \citep{schick2023toolformer}. A recent transport-planning review documented the expanding role of large language models in modelling, design and decision making \citep{jin2026a}. These developments support a conversational interface, but they do not remove the need for deterministic analysis.

ResiliFlow treats the language model as an orchestrator. It can interpret a request, recommend a tool and explain a recorded result. Network metrics, routes, assignments, simulations and detector outputs remain products of executable modules. Reproducibility rules record the input, model version, parameters, checks and artefacts for each run \citep{pritchard2025formal}. A local Assistant provides guidance without an external key. An optional Copilot supports multiple providers when a user supplies a valid key. This separation keeps the platform usable when a provider is unavailable and prevents fluent text from becoming unrecorded model evidence.

\subsection{The integration gap}

Existing approaches usually optimise one part of the workflow. Disaster dashboards display hazard and network layers. Digital twins synchronise data and assets. AI inspection platforms detect visible damage. Simulation systems test behaviour under specified scenarios. Conversational assistants interpret questions. ResiliFlow connects these capabilities through one shared state and makes their boundaries explicit. Table~\ref{tab:positioning} summarises the intended distinction.

\begingroup
\footnotesize
\setlength{\tabcolsep}{3pt}

\begin{longtable}{@{}L{2.6cm}L{3.0cm}L{5.0cm}L{4.2cm}@{}}
\caption{Capability positioning of ResiliFlow.}
\label{tab:positioning}\\
\toprule
Approach
& Primary strength
& Connection to state and decisions
& Typical endpoint \\
\midrule
\endfirsthead

\toprule
Approach
& Primary strength
& Connection to state and decisions
& Typical endpoint \\
\midrule
\endhead

\bottomrule
\endfoot

Disaster dashboard
& Fast situation display
& Often presents selected layers and indicators
& Operational picture \\

Urban digital twin
& Synchronised asset or city representation
& Maintains a live or frequently updated digital state
& Asset or city view \\

AI road inspection
& Scalable visual condition detection
& Converts imagery into condition observations
& Defect detections \\

Transport simulator
& Behaviour under specified scenarios
& Converts assumptions into trajectories and performance measures
& Scenario comparison \\

Conversational assistant
& Natural-language access
& Interprets questions and explains connected information
& Guided interaction \\

\textbf{ResiliFlow}
& \textbf{Connected observation, modelling and action}
& \textbf{Maintains state, executes models, records evidence and supports feedback}
& \textbf{Reviewable decision evidence} \\

\end{longtable}
\endgroup
\FloatBarrier

\section{The ResiliFlow transport world model platform}

\subsection{An operational definition}

\begin{figure}[H]
  \centering
  \includegraphics[width=0.98\textwidth]{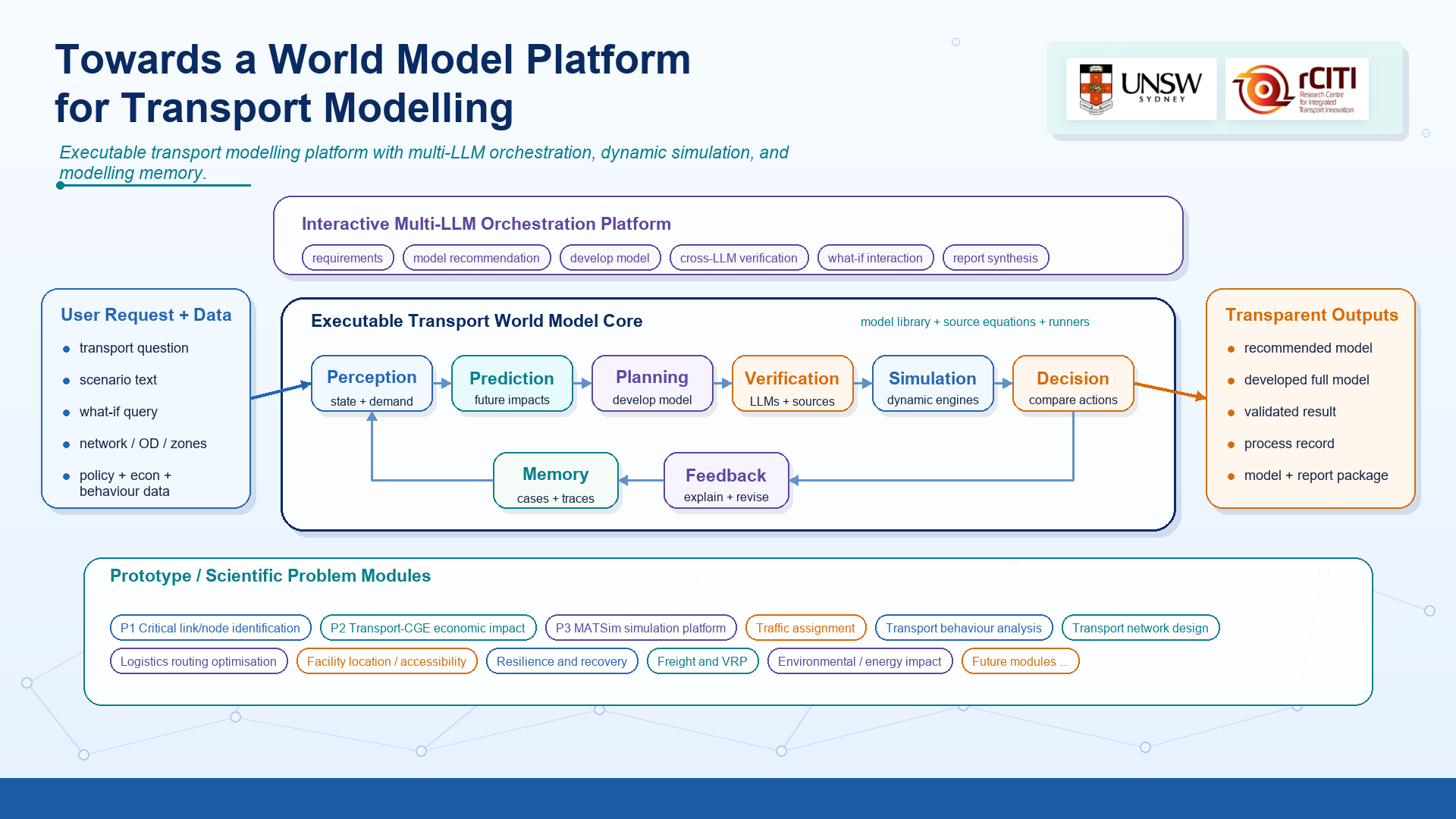}
  \caption{ResiliFlow concept architecture. The diagram presents the long-term transport world model structure and the relationship between user requests, executable modelling, verification, outputs and reusable memory.}
  \label{fig:concept}
\end{figure}

Figure~\ref{fig:concept} presents the macro architecture supplied during the platform design process. It shows how user requests and data enter an executable core, how the model library supports the core and how transparent outputs leave it. The feedback and memory links make the cycle reusable.

\FloatBarrier

\subsection{Eight connected steps}

\begin{figure}[H]
  \centering
  \includegraphics[width=0.98\textwidth]{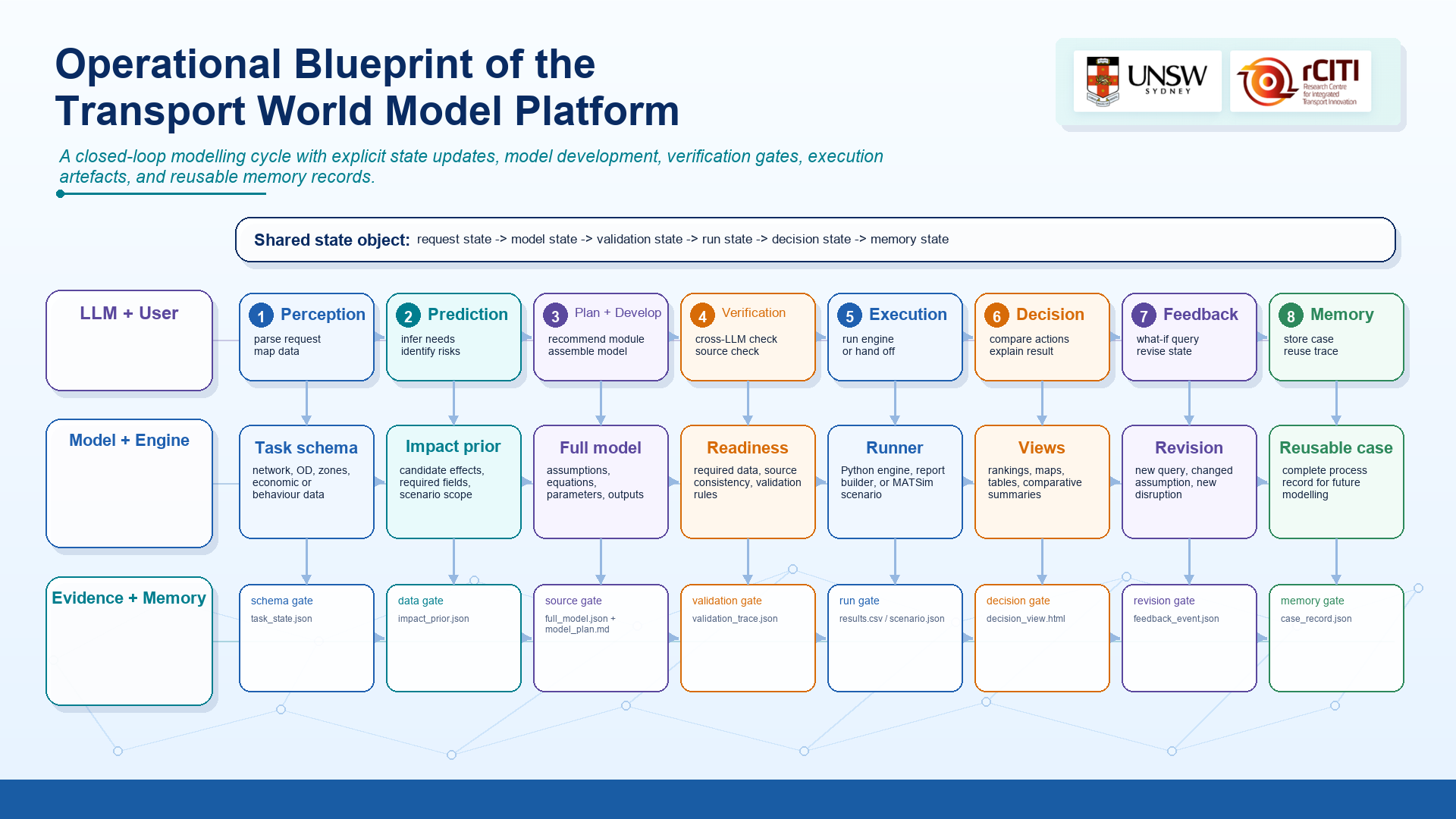}
  \caption{Operational blueprint of the transport world model. The shared state links eight stages to their executable artefacts and evidence gates.}
  \label{fig:blueprint}
\end{figure}

The operational cycle is shown in Figure~\ref{fig:blueprint}. Perception maps available data into a current state. Prediction identifies likely needs and consequences. Planning assembles a model and its assumptions. Verification checks source, schema and readiness. Execution runs an engine or records a human hand-off. Decision turns results into comparable views. Feedback changes the question, data or assumptions. Memory stores a reusable case.

\FloatBarrier

This cycle allows different users to enter at different points. An emergency planner may begin with a disruption question. An engineer may begin with an image. A researcher may begin with a model specification. The same state structure records what happened next. This is the main practical difference between a workflow diagram and a world model platform. The steps are executable and each material transition leaves evidence.

\subsection{Shared state, model library and evidence gates}

The shared state has four responsibilities. It identifies the active project, stores mapped entities and time coverage, records selected assumptions and links outputs to their generating run. The model library provides transparent equations, parameter definitions and runners. Evidence gates check whether the required fields, sources and validation conditions are present. Memory stores a complete case only when the user chooses to save it.

The platform distinguishes three output classes. \emph{Descriptive evidence} reports what is present in the data. \emph{Analytical evidence} is produced by a documented model. \emph{Decision evidence} adds priorities, trade-offs and human review. This separation prevents a map layer, a model prediction and a recommendation from appearing equivalent.

\subsection{Question-led orchestration}

The local Assistant explains tools and data requirements. The optional LLM Copilot accepts open questions such as ``Which road failures most reduce access?'' and ``Which repair sequence restores the greatest network performance?'' It converts the request into a task schema, recommends one or more modules and asks for missing data. It does not calculate network metrics in free text. The executable module produces the result, and the Copilot explains that result with links to the recorded evidence.

This architecture supports both single and compound tool use. A routing request can call one function. A recovery question can combine critical-link analysis, area aggregation, repair sequencing and a scenario comparison. The eight-step Research Validation view remains available as a fixed scientific record. Future question-led validation will expose the same stages in a form tailored to a user's research question.

\section{Disaster Transport Resilience Analysis}

\subsection{A data-first workspace}

The first workspace begins with user data. A user can upload one compatible file or combine several files. The platform maps the records before it activates an analysis. It accepts several evidence tiers, ranging from a ready network to observations that must be converted into links, nodes, zones or origin--destination movements. OpenStreetMap can provide an openly documented street-network starting point when its provenance and fitness for purpose are recorded \citep{haklay2008openstreetmap}. Pervasive traffic observations can support origin--destination estimation and related inputs \citep{waller2021rapidex}. No demonstration project is opened automatically. Temporary projects remain local until the user explicitly saves them.

\begin{figure}[H]
  \centering
  \includegraphics[width=0.98\textwidth]{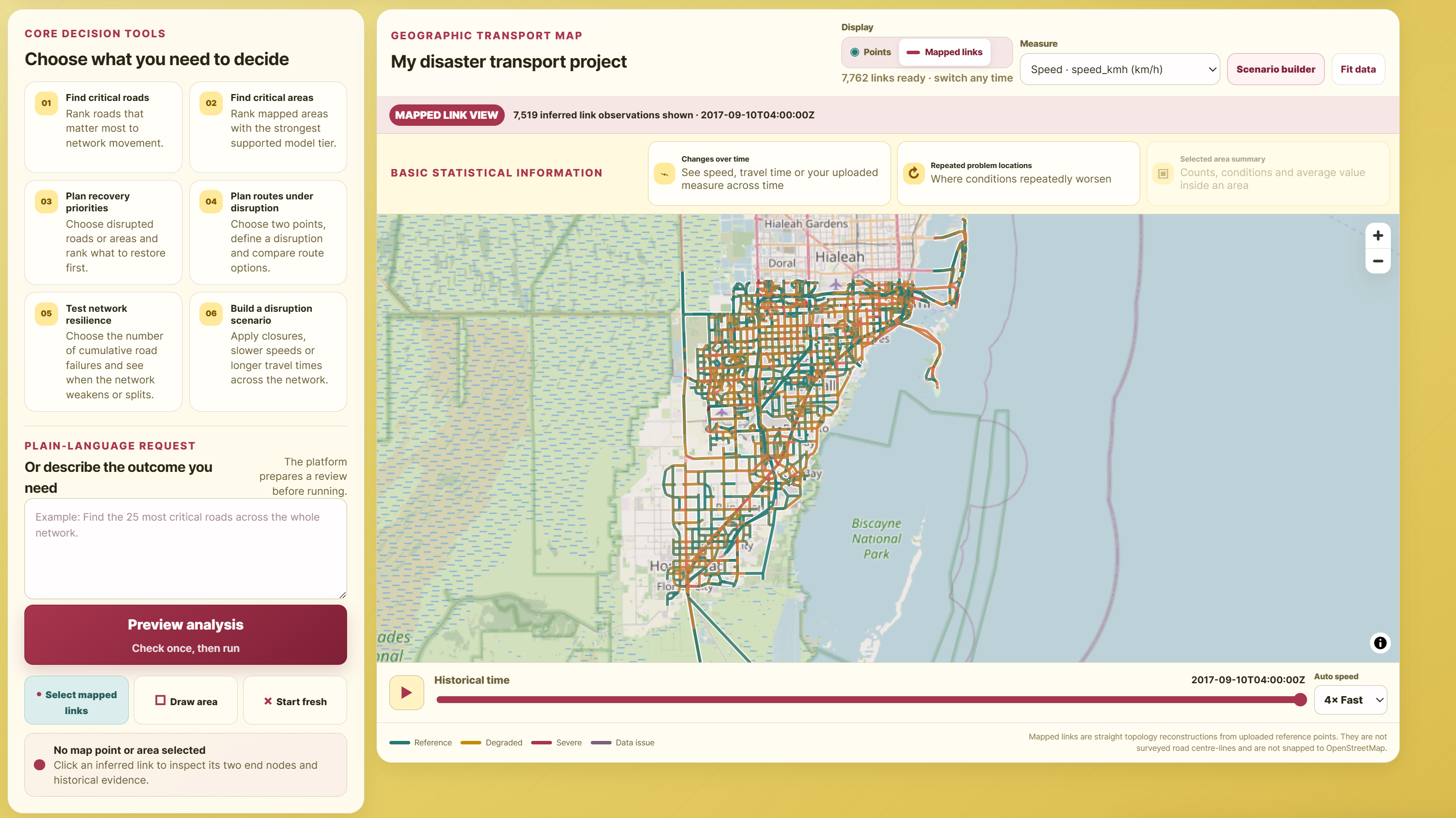}
  \caption{The implemented Disaster Transport Resilience Analysis workspace. Six decision tools share the same geographical map, active data, historical-time control and project state.}
  \label{fig:disasterworkspace}
\end{figure}

The current implementation contains six core tools within one map-centred workspace. Users can select a tool directly or describe the outcome they need. The map, time control, active measure and project state remain visible while they move between critical-road analysis, critical-area analysis, recovery planning, disruption routing, resilience testing and scenario simulation. Figure~\ref{fig:disasterworkspace} shows this implemented workspace with the six decision tools and a mapped transport network.

Table~\ref{tab:six-tools} connects each operational question to representative model logic and the output that a user sees. The shared workspace allows the result of one function to become the evidence or selection context for another.

\begin{table}[H]
\centering
\caption{Six core analysis functions, representative model logic and user-facing outputs.}
\label{tab:six-tools}
\footnotesize
\begin{tabularx}{\textwidth}{L{3.0cm}L{4.3cm}Y}
\toprule
Function & Representative model logic & User-facing output \\
\midrule
Critical roads & Link-removal change in network efficiency, connectivity or service & Ranked links and mapped criticality \\
Critical areas & Spatial aggregation of link or node evidence & Hotspots with contributing records \\
Recovery priorities & Marginal performance gain per repair cost or resource unit & Ordered repair sequence and benefit trace \\
Disruption routing & Shortest and alternative paths under closures, delay and risk & Viable routes, costs and affected links \\
Resilience testing & Targeted or random removal with retained-performance tracking & Resilience curve and comparative score \\
Scenario simulation & Closures, speed loss, assignment and optional agent-based execution & Before-and-after scenario outcomes \\
\bottomrule
\end{tabularx}
\end{table}

\subsection{Development of the Model Library}
\label{sec:model-library}

ResiliFlow includes an extensible model library that connects user questions, available data and suitable analytical methods. The library supports the transition from transport observations to executable analysis. It also allows different analytical functions to share a common project state. This design helps users move between geographic summaries, network diagnosis, disruption analysis, routing, recovery planning and simulation without rebuilding the project for every task. The current model library supports geographic condition summaries, recurring-problem detection, critical-road identification, selected-road closure analysis, recovery prioritisation, disruption scenario construction, disruption propagation and routing under disrupted conditions. These models provide the analytical foundation for the user-facing tools in the Disaster Transport Resilience Analysis workspace. They also support the LLM Copilot by providing executable methods that can be selected and run in response to plain-language questions.

\begin{figure}[H]
  \centering
  \includegraphics[width=0.98\textwidth]{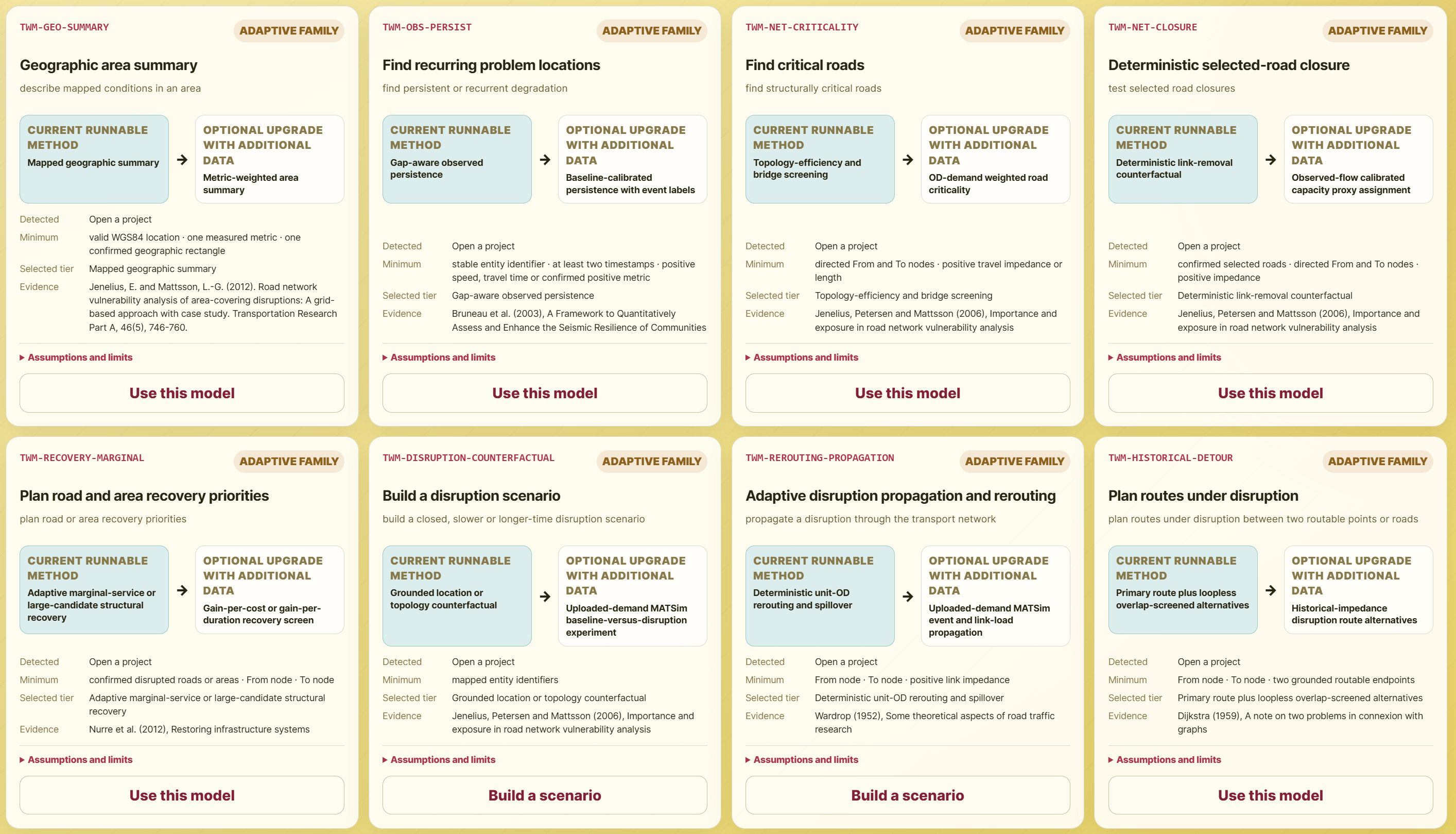}
  \caption{Representative excerpt from the ResiliFlow model library. The illustrated models support geographic condition summaries, recurring-problem detection, critical-road identification, selected-road closure analysis, recovery prioritisation, disruption scenario construction, disruption propagation and routing under disrupted conditions. Each model card links an operational method with its data requirements, potential extensions and supporting scientific evidence.}
  \label{fig:model-library}
\end{figure}

Figure~\ref{fig:model-library} presents a representative excerpt from the ResiliFlow model library. It illustrates how each model card connects a transport question with an operational method, required data, potential extensions and supporting scientific evidence.

\FloatBarrier

\subsection{Critical Roads and Areas}

\begin{figure}[H]
  \centering
  \includegraphics[width=0.97\textwidth]{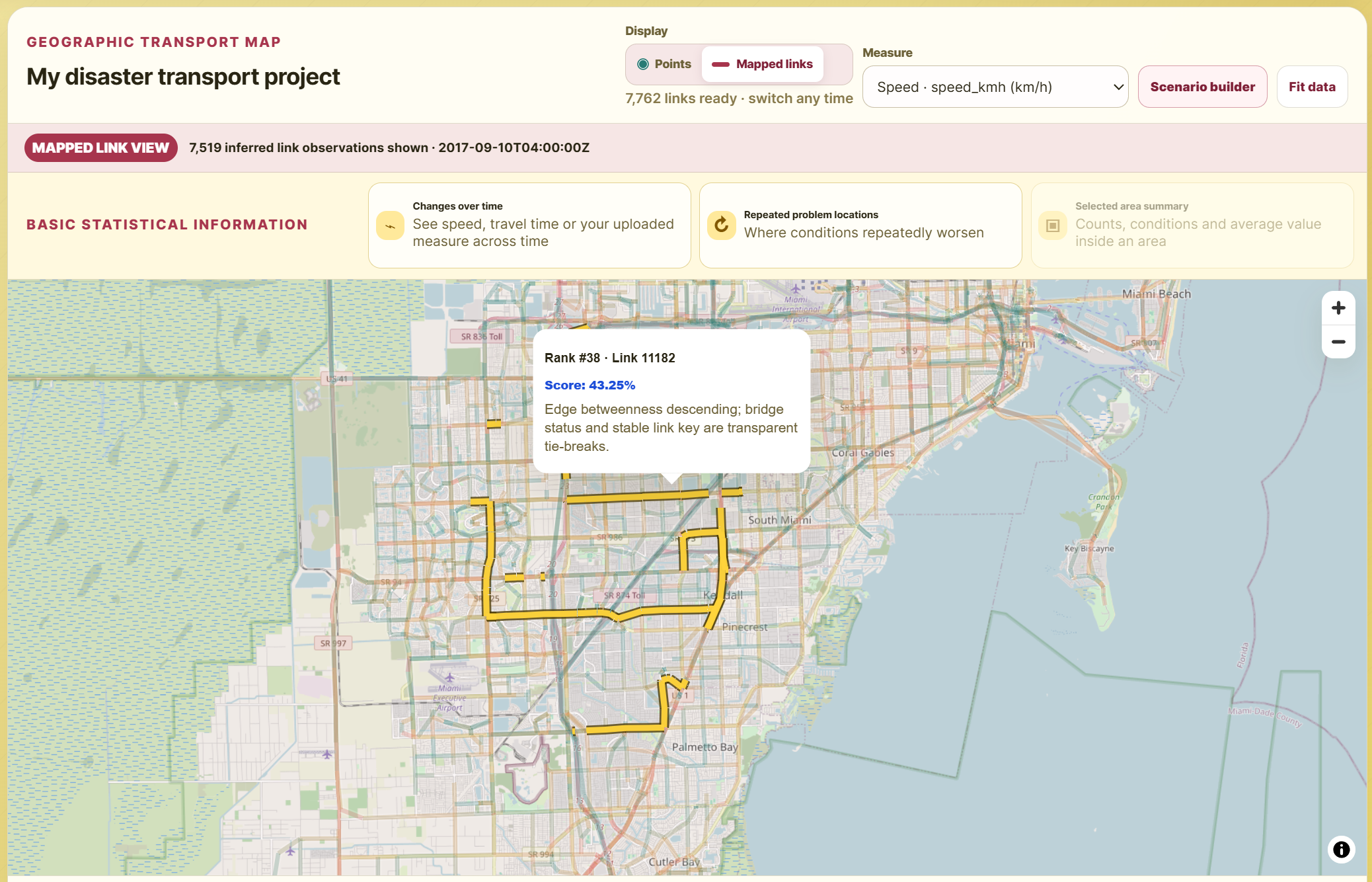}
  \caption{Implemented critical-link result in the shared geographical workspace.}
  \label{fig:criticallink}
\end{figure}

\paragraph{Critical roads.}
Critical-road analysis identifies and ranks roads whose disruption may cause the greatest loss of network connectivity or performance. As shown in Figure~\ref{fig:criticallink}, ResiliFlow selects a suitable method according to the available network, travel and demand data.

\FloatBarrier

\begin{figure}[H]
  \centering
  \includegraphics[width=0.97\textwidth]{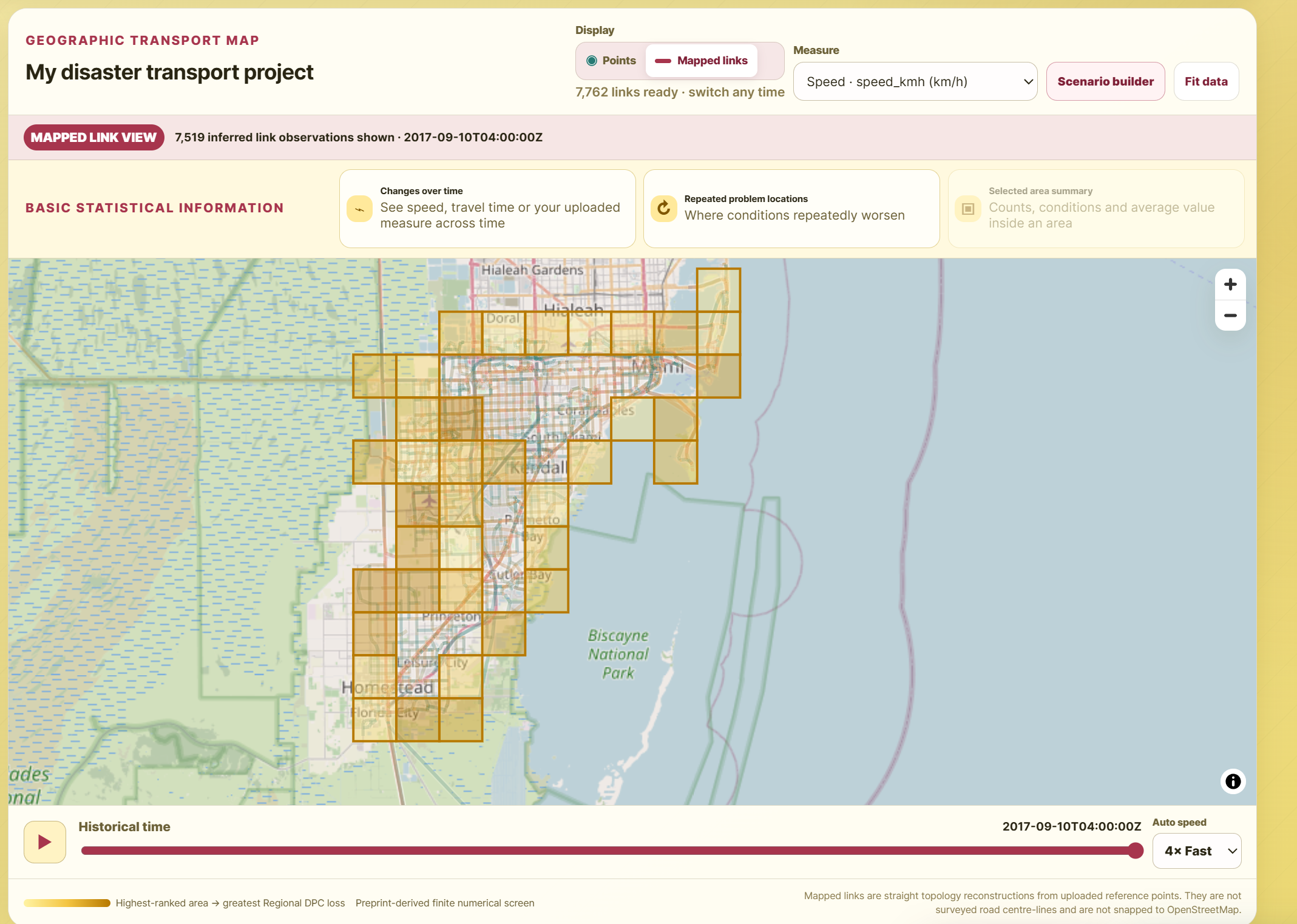}
  \caption{Implemented critical-area result in the shared geographical workspace.}
  \label{fig:criticalarea}
\end{figure}

\paragraph{Critical areas.}
Critical-area analysis aggregates road and node information within mapped spatial units to locate areas with high disruption exposure or transport importance. As shown in Figure~\ref{fig:criticalarea}, the platform can summarise geographic conditions and support more detailed weighting when additional data become available.

\FloatBarrier

\subsection{Recovery Priorities}

\begin{figure}[H]
  \centering
  \includegraphics[width=0.97\textwidth]{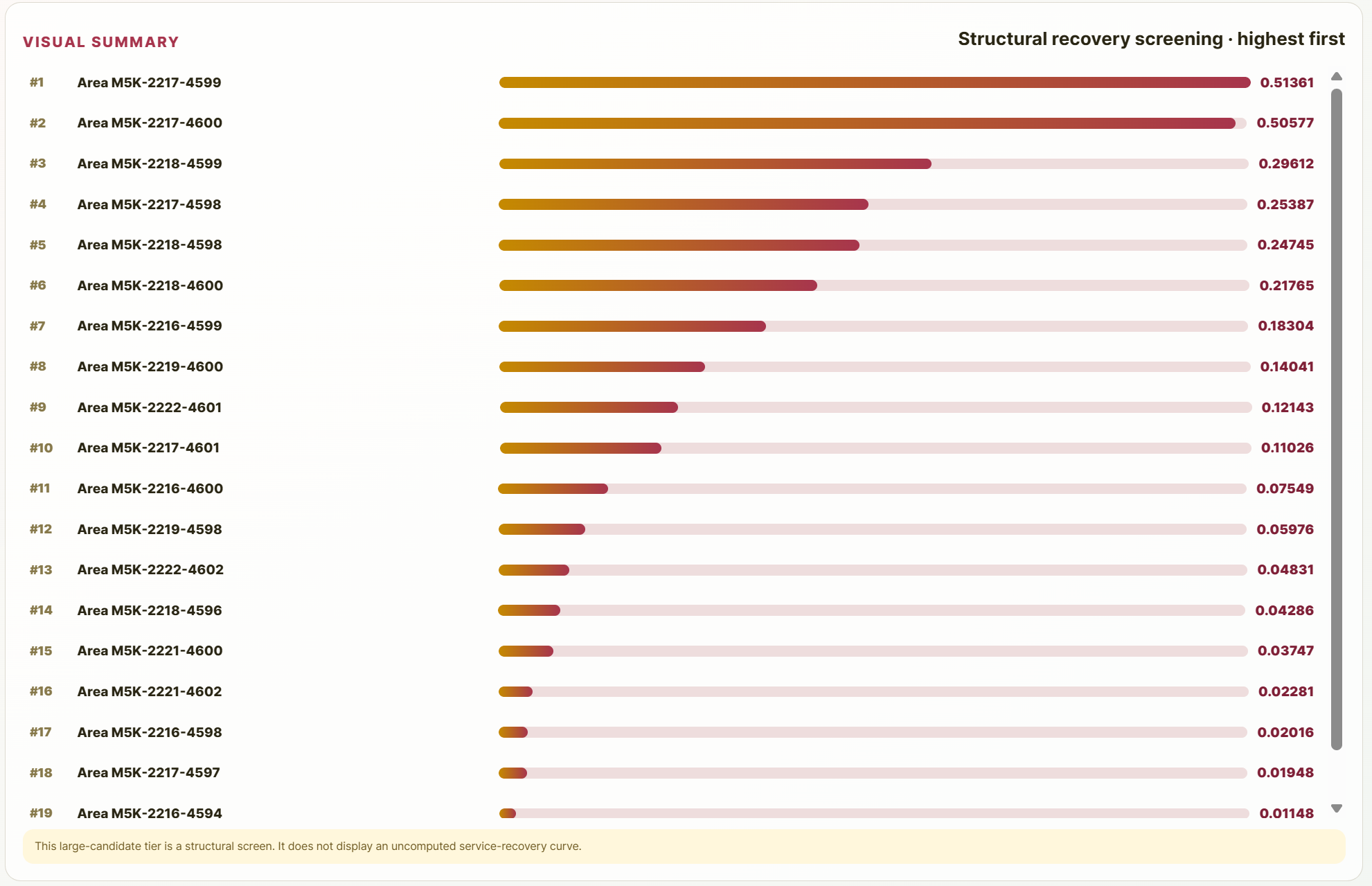}
  \caption{Implemented critical-area recovery prioritisation result in the shared geographical workspace.}
  \label{fig:recoveryprio}
\end{figure}

Recovery-priority analysis ranks damaged roads or areas according to the expected benefit of restoration under limited budgets and resources. The ranking can consider network performance, repair cost, accessibility, reliability and equity \citep{xu2024exploring,xu2024predisaster,xu2025predisaster,niu2023predisaster,zhang2023equity,najmi2023equity,zhang2025integrating}. As shown in Figure~\ref{fig:recoveryprio}, ResiliFlow provides a runnable recovery method and supports methodological upgrades when richer cost, demand or service data become available.

\FloatBarrier

\subsection{Disruption Routing}

\begin{figure}[H]
  \centering
  \includegraphics[width=0.97\textwidth]{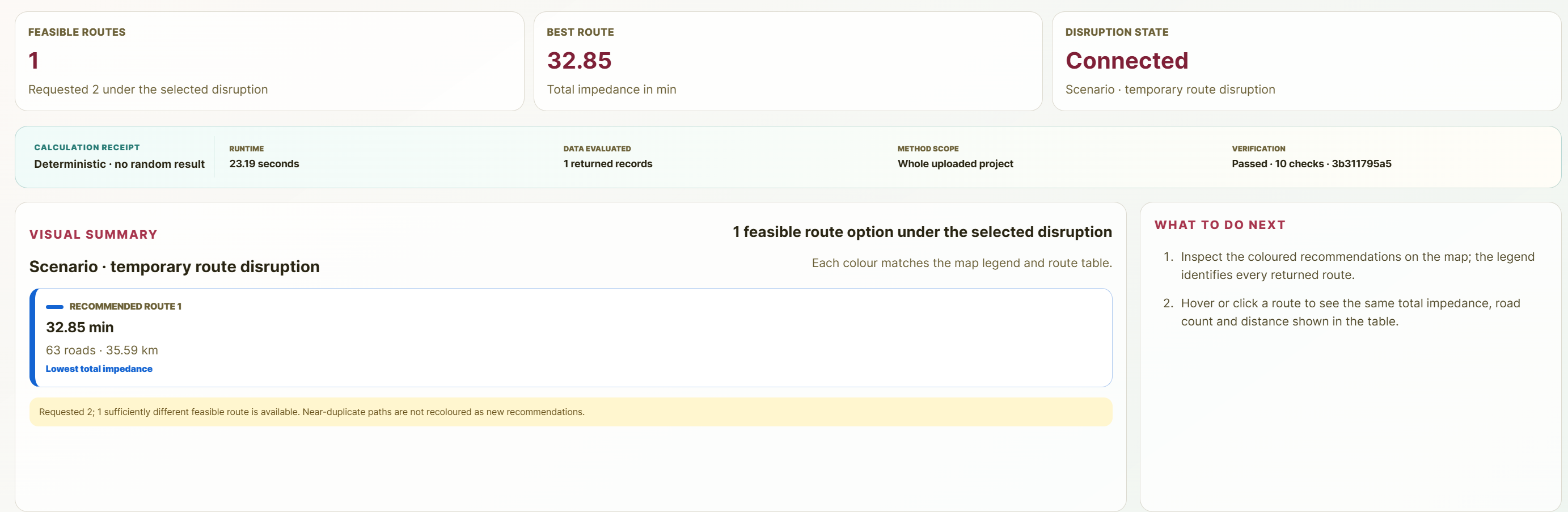}
  \caption{Optimal route recommendations under different road scenarios.}
  \label{fig:disruroutes}
\end{figure}

Disruption-routing analysis identifies viable paths after road closures or changes in travel conditions. ResiliFlow can generate a primary route and alternative paths that avoid disrupted or exposed links using established shortest-path methods \citep{dijkstra1959a,yen1971finding}. As shown in Figure~\ref{fig:disruroutes}, the routing models can incorporate richer historical and disruption information when suitable data are available.

\FloatBarrier

\subsection{Resilience testing}

Resilience testing removes links or nodes according to a chosen strategy and traces retained performance. Let $G_k$ be the network after $k$ removals. A normalised curve is
\begin{equation}
R_k=\frac{E(G_k)}{E(G_0)}.
\label{eq:resilience}
\end{equation}
The area under this curve summarises how quickly performance declines. Targeted and random removals answer different questions. A targeted test examines vulnerability to informed disruption. A random test provides a baseline. Real-world disruption studies have shown why scenario selection, limited information and observed network conditions matter \citep{bagloee2017identifying,waller2025rapid,niu2022linklevel}.

\begin{figure}[H]
  \centering
  \includegraphics[width=0.97\textwidth]{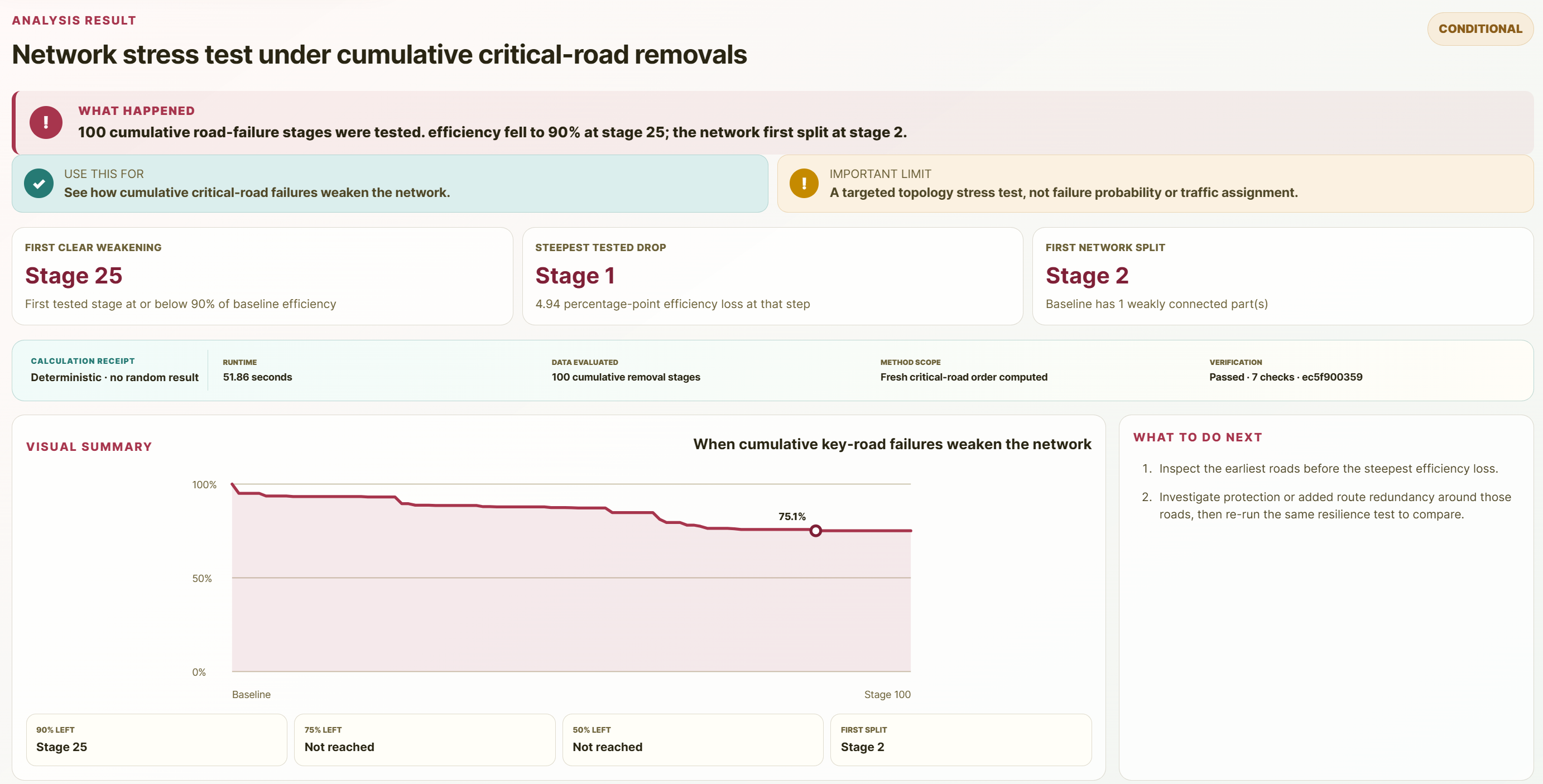}
  \caption{Resilience testing under cumulative road disruptions in ResiliFlow.}
  \label{fig:resilience}
\end{figure}

As shown in Figure~\ref{fig:resilience}, the resulting resilience curve supports comparison of disruption scenarios and recovery needs.

\FloatBarrier

\subsection{Scenario simulation}

The scenario builder applies closures, speed reductions, travel-time increases and demand changes. A speed-loss rate $r_e$ changes free-flow time to
\begin{equation}
t'_e=\frac{t_e}{1-r_e},
\label{eq:speedloss}
\end{equation}
while a direct time multiplier $m_e$ gives $t'_e=t_e(1+m_e)$. Static assignment can use the Bureau of Public Roads function
\begin{equation}
t_e(x_e)=t_e^0\left[1+\alpha\left(\frac{x_e}{c_e}\right)^{\beta}\right],
\label{eq:bpr}
\end{equation}
with all-or-nothing assignment or Frank--Wolfe iterations \citep{beckmann1956studies,bureau1964traffic,frank1956an}. MATSim provides an optional agent-based execution path for scenarios that need population-level activity and route adaptation \citep{axhausen2016the}. Evacuation dynamics, regional validation and simulation studies guide the interpretation of such runs \citep{dixit2011validation,dixit2014evacuation,lalwani2026traffic,sevim2025a}.

Broader economic effects can be handled by linked models. Integrated transport and computable general equilibrium research has established methods for connecting network change to economic activity \citep{robson2017a,robson2018a,shahriari2023integrating,wang2024calibration}. Related equilibrium research has also connected transport technologies, parking behaviour and economic response \citep{zhang2021an}. ResiliFlow treats these as model-library extensions with their own data and validation needs.

\begin{figure}[H]
  \centering
  \includegraphics[width=0.97\textwidth]{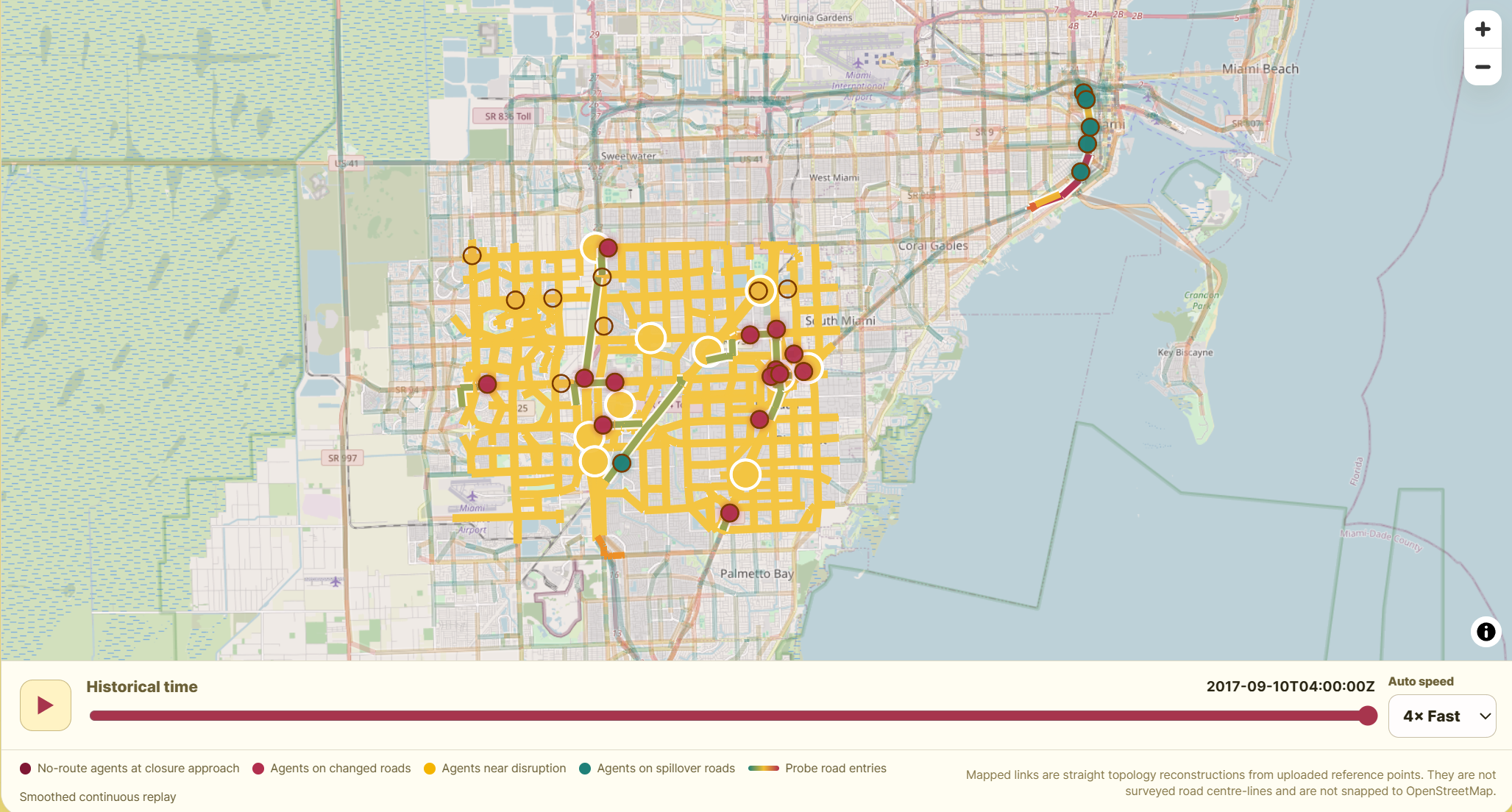}
  \caption{The MATSim simulation process in ResiliFlow.}
  \label{fig:matsim2}
\end{figure}

The MATSim simulation process for the specified scenario is shown in Figure~\ref{fig:matsim2}.

\subsection{Research Validation, memory and transparent outputs}

Each tool reports more than a headline number. The Research Validation module records the question, data readiness, assumptions, source equations, execution status and limitations. A result view can contain a map, ranking, comparison chart, table and downloadable run record. The user can revise the question or assumptions and retain the relationship between versions.

\begin{figure}[t]
  \centering
  \includegraphics[width=0.97\textwidth]{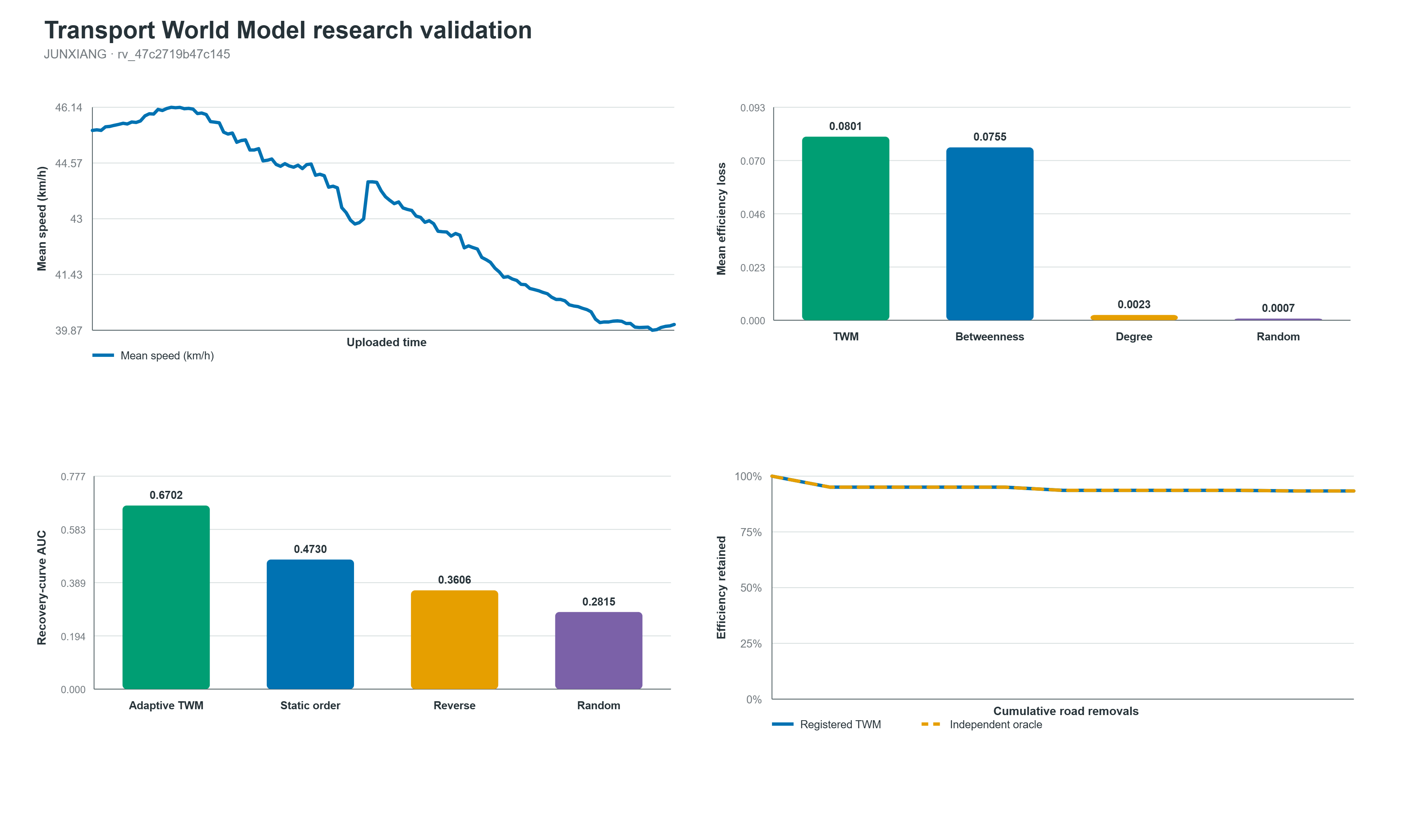}
  \caption{Research Validation output. The platform presents time-series evidence, benchmark comparison, recovery-curve performance and retained network efficiency in one visual record.}
  \label{fig:validationresults}
\end{figure}

\begin{figure}[t]
  \centering
  \includegraphics[width=0.97\textwidth]{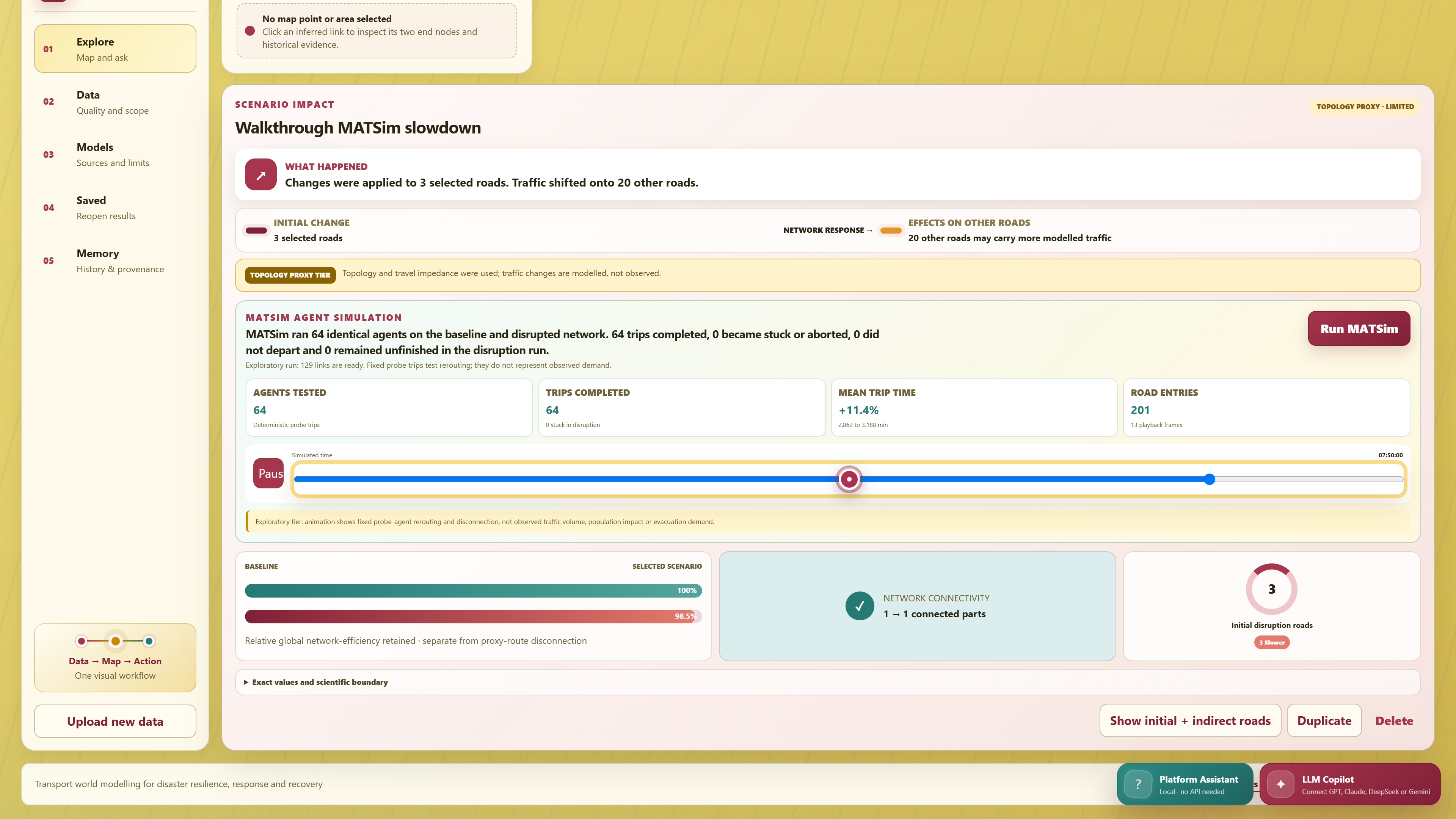}
  \caption{MATSim scenario result embedded in the disaster-resilience workspace. The interface links the initial road change to agent outcomes, network response and interpretable performance indicators.}
  \label{fig:matsim}
\end{figure}

Figure~\ref{fig:validationresults} illustrates how validation turns model outputs into readable comparative evidence. The four views report observed speed, efficiency loss, recovery performance and retained network efficiency. Figure~\ref{fig:matsim} shows the optional MATSim pathway within the same workspace. It connects the selected disruption to agent trips, travel-time change, network connectivity and indirect road effects. 

\section{AI-based Transport Infrastructure Perception and Decision Support}

\subsection{From imagery to usable condition evidence}

The second workspace begins at a location. A user can inspect a road through interactive street-level imagery, search historical satellite context and add authorised ground evidence. The interface keeps Google imagery for interactive human review. It does not use ordinary Google Maps or Street View access as an AI training licence. KartaView and audited public datasets are candidates only after source, attribution, privacy and licence checks. Copernicus Sentinel imagery supports a separate hazard-context function because its resolution is not suitable for street-level crack detection.

\begin{figure}[H]
  \centering
  \includegraphics[width=0.98\textwidth]{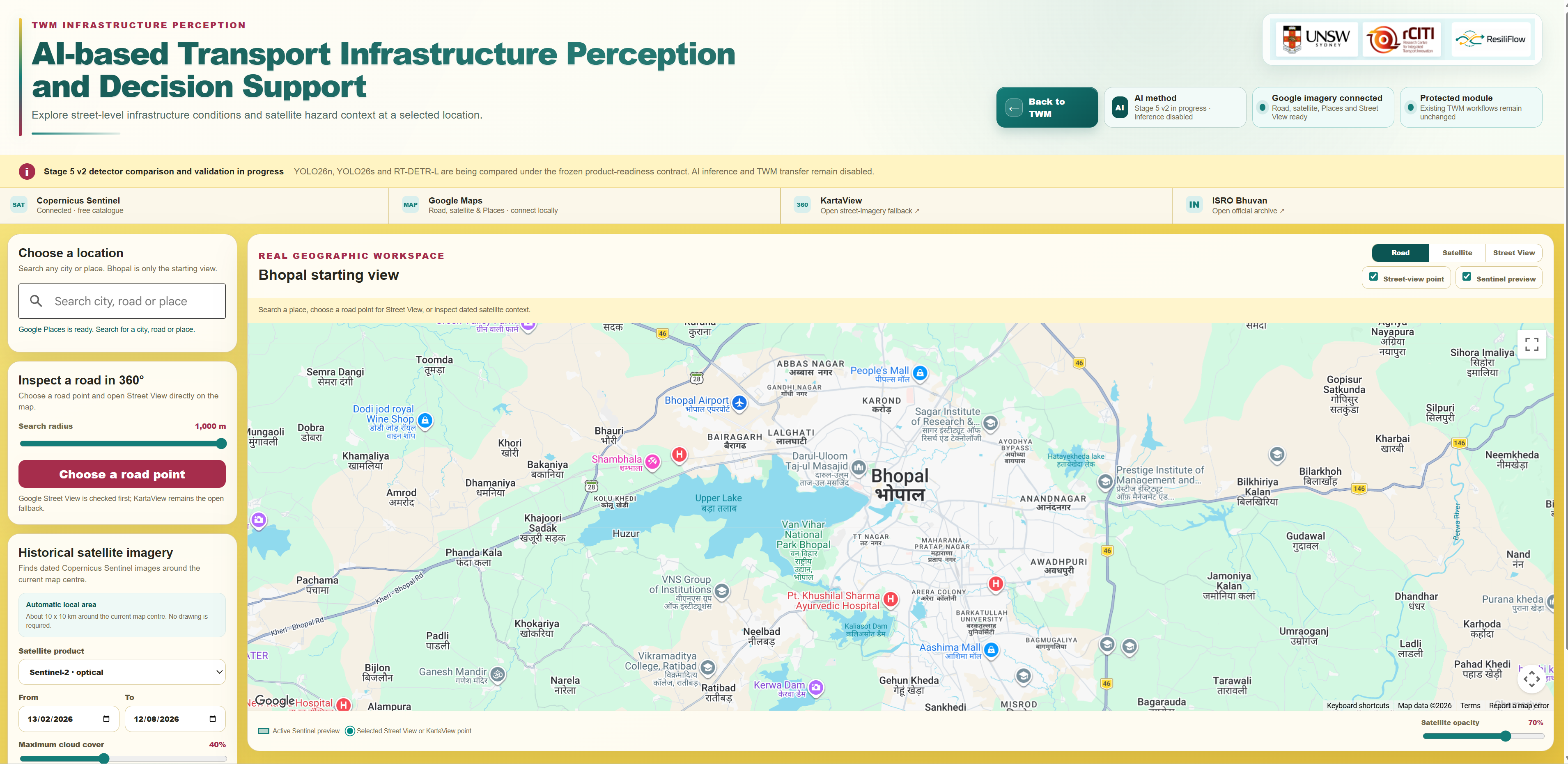}
  \caption{Complete AI-based Transport Infrastructure Perception and Decision Support workspace centred on Bhopal. Location search, 360-degree inspection, historical satellite imagery and the geographic workspace remain visible together.}
  \label{fig:perceptionui}
\end{figure}

Figure~\ref{fig:perceptionui} shows the complete workspace. The user can search a city, select a road, move from map context into a street view and call historical satellite context without leaving the location. This geographic continuity gives the detector and the decision layer a shared place-based state.

\begin{figure}[H]
  \centering
  \includegraphics[width=0.98\textwidth]{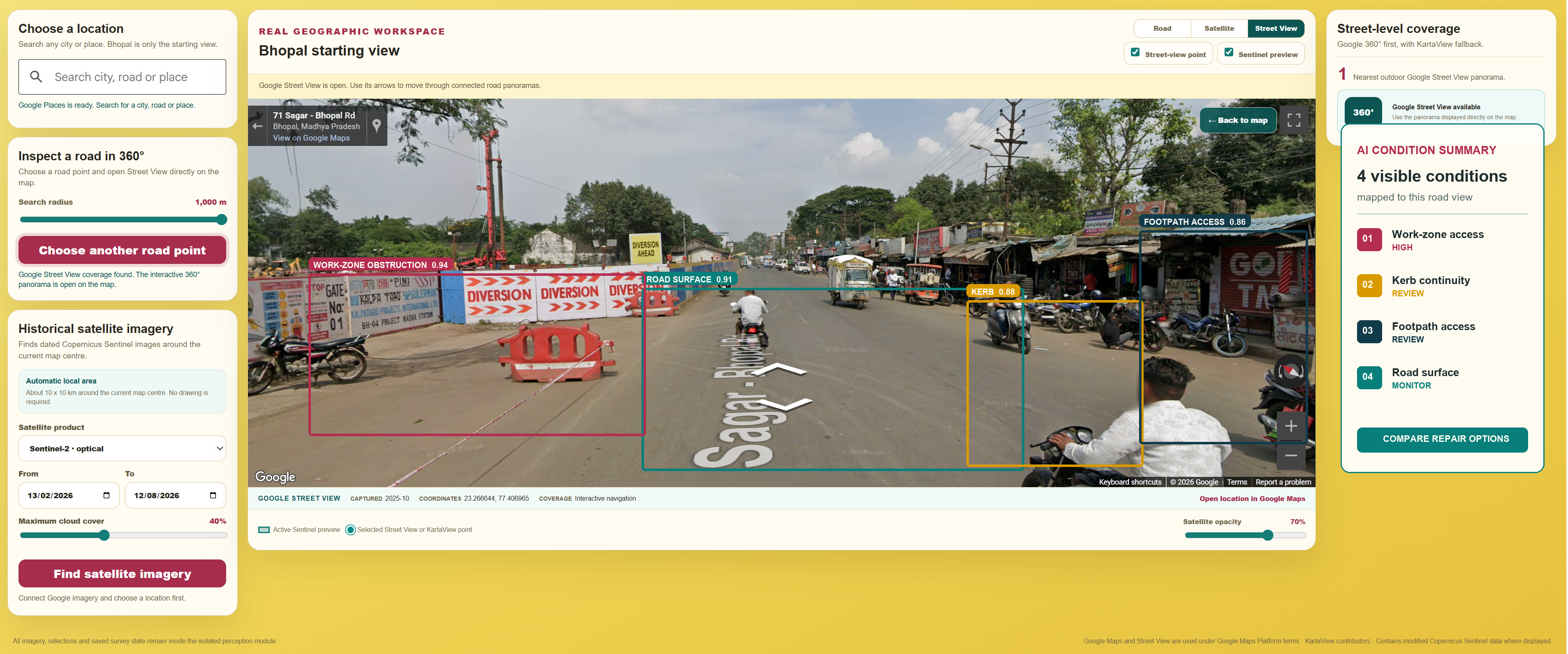}
  \caption{Bhopal street-view perception and decision-support experience. AI observations appear in the geographic workspace and lead directly to a concise condition summary and repair-option comparison. Google attribution and the interactive street-view context are retained.}
  \label{fig:bhopalai}
\end{figure}

Figure~\ref{fig:bhopalai} presents the intended experience after a road view is selected. The visual perception layer marks relevant road, footpath, kerb and access conditions. The decision card groups the observations into actions that can be reviewed, compared and passed to repair planning. The overlay is a platform interface visualisation. It demonstrates how detector outputs are designed to appear in the workspace.

\begin{figure}[H]
  \centering
  \includegraphics[width=0.92\textwidth]{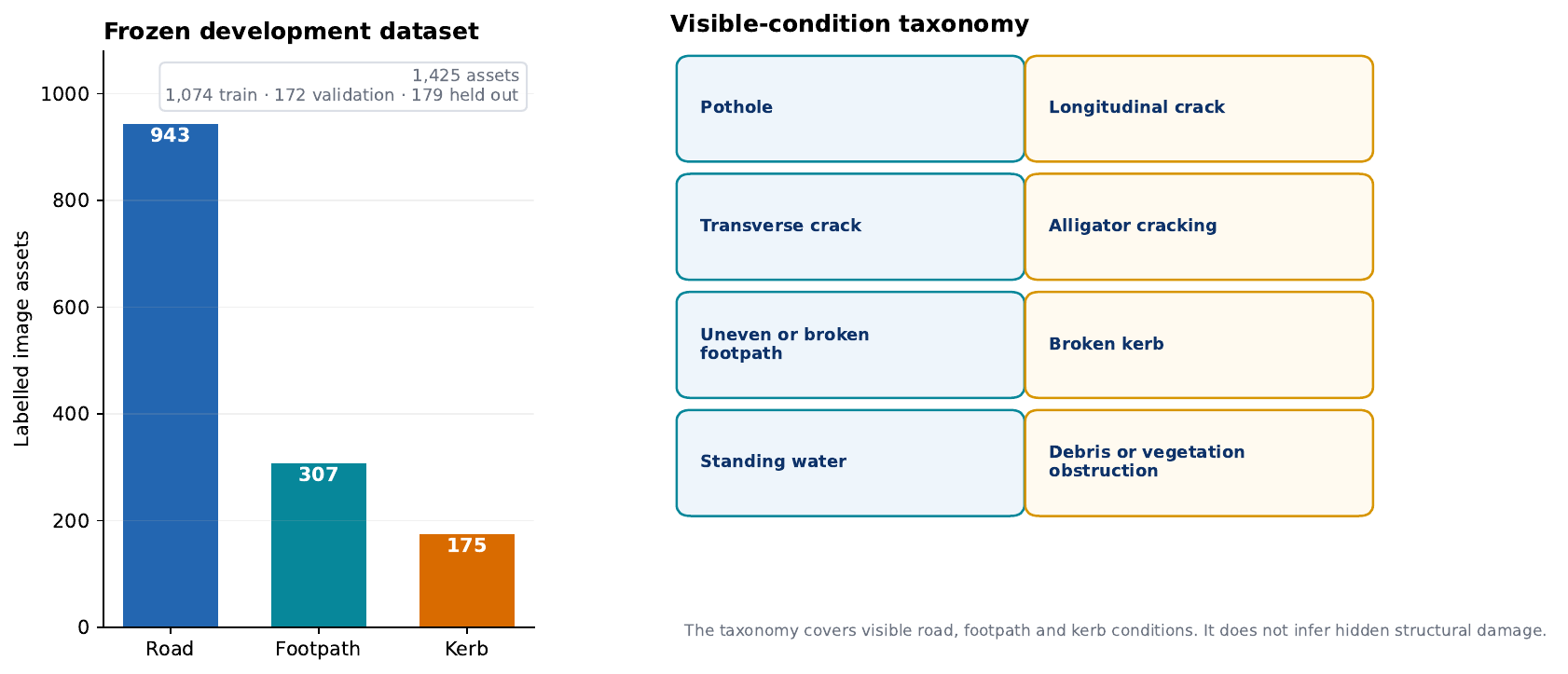}
  \caption{Frozen baseline development dataset and the visible-condition taxonomy. The 1,425 labelled assets comprise 1,074 training, 172 validation and 179 held-out images across road, footpath and kerb evidence.}
  \label{fig:dataset}
\end{figure}

The first condition taxonomy contains eight visible classes: pothole, longitudinal crack, transverse crack, alligator cracking, uneven or broken footpath, broken kerb, standing water, and debris or vegetation obstruction. Table~\ref{tab:conditions} translates these classes into practical inspection and decision uses. Figure~\ref{fig:dataset} summarises the frozen baseline dataset and taxonomy. The labels cover visible surface conditions. They do not claim hidden structural damage.

\begingroup
\footnotesize
\setlength{\tabcolsep}{3pt}

\begin{longtable}{@{}L{3.5cm}L{3.0cm}L{8.4cm}@{}}
\caption{Visible-condition classes and how they can inform infrastructure decisions.}
\label{tab:conditions}\\
\toprule
Visible condition
& Infrastructure element
& Illustrative decision use \\
\midrule
\endfirsthead

\toprule
Visible condition
& Infrastructure element
& Illustrative decision use \\
\midrule
\endhead

\bottomrule
\endfoot

Pothole
& Road
& Prioritise safety inspection and patching \\

Longitudinal crack
& Road or footpath
& Track linear deterioration and maintenance need \\

Transverse crack
& Road or footpath
& Review cross-surface failure and ride or access quality \\

Alligator cracking
& Road
& Identify concentrated pavement fatigue \\

Uneven or broken footpath
& Footpath
& Screen accessibility and trip-hazard concerns \\

Broken kerb
& Kerb
& Review edge continuity, drainage and pedestrian access \\

Standing water
& Road, footpath or kerb
& Connect visible ponding to drainage and hazard review \\

Debris or vegetation obstruction
& Road, footpath or kerb
& Identify clearance and continuity priorities \\

\end{longtable}
\endgroup
\FloatBarrier

\subsection{AI detection and learning}

The first reproducible model was a compact convolutional detector with 1,258,087 parameters and a $160\times160$ input. It was trained under three seeds on a frozen split. The held-out test contained 179 images. Its macro F1 was 0.573, mean intersection over union was 0.491, mAP@0.5 was 0.268, mAP@0.5:0.95 was 0.151 and ten-bin expected calibration error was 0.116. This baseline established the full path from lawful training data to localised condition evidence.

The next dataset contains 15,596 unique training images and 2,644 validation images, with weighted training exposure equivalent to 21,877 samples. It combines global road evidence, Indian road evidence, a smaller Indian footpath set, secondary sources and background images. Three candidate families were selected for full comparison: You Only Look Once (YOLO)26n, YOLO26s and Real-Time Detection Transformer, Large variant (RT-DETR-L). YOLO26n and YOLO26s test the speed--capacity trade-off within a modern real-time detector family \citep{jocher2026yolo26}. RT-DETR-L provides a transformer-based end-to-end alternative \citep{zhao2024detrs}.

\begin{figure}[H]
  \centering
  \includegraphics[width=0.91\textwidth]{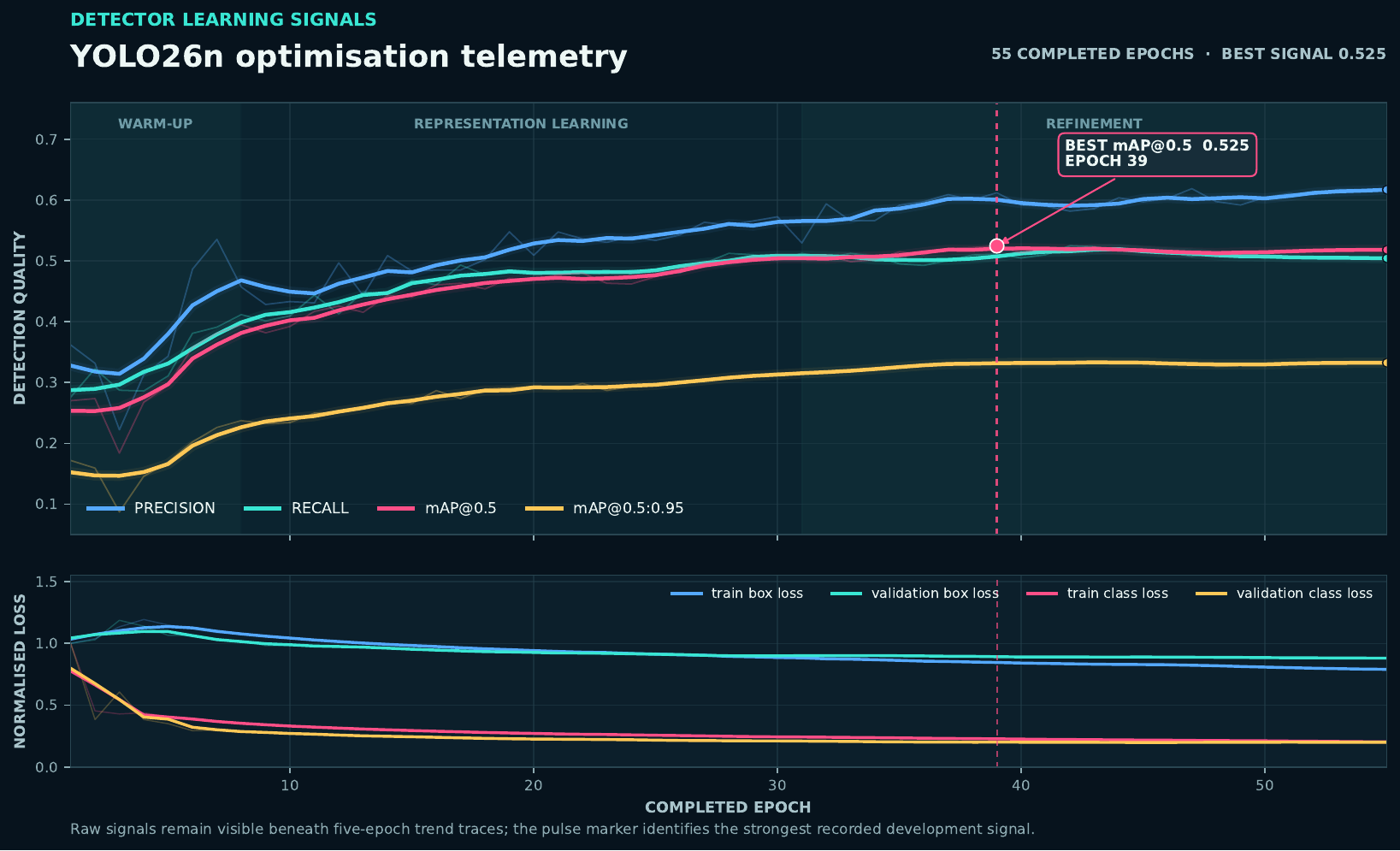}
  \caption{YOLO26n learning curve across precision, recall and two average-precision measures. The progression shows how the detector develops a useful balance between finding conditions and limiting false detections.}
  \label{fig:training}
\end{figure}

Figure~\ref{fig:training} records the available YOLO26n run. The strongest recorded development-set mAP@0.5 was 0.525. The wider comparison uses YOLO26n and YOLO26s for efficient real-time detection and RT-DETR-L as a transformer alternative. Together they provide a practical route to both responsive platform use and higher-capacity perception.

\subsection{Recognition results and cross-source learning}

\begin{figure}[t]
  \centering
  \includegraphics[width=0.90\textwidth]{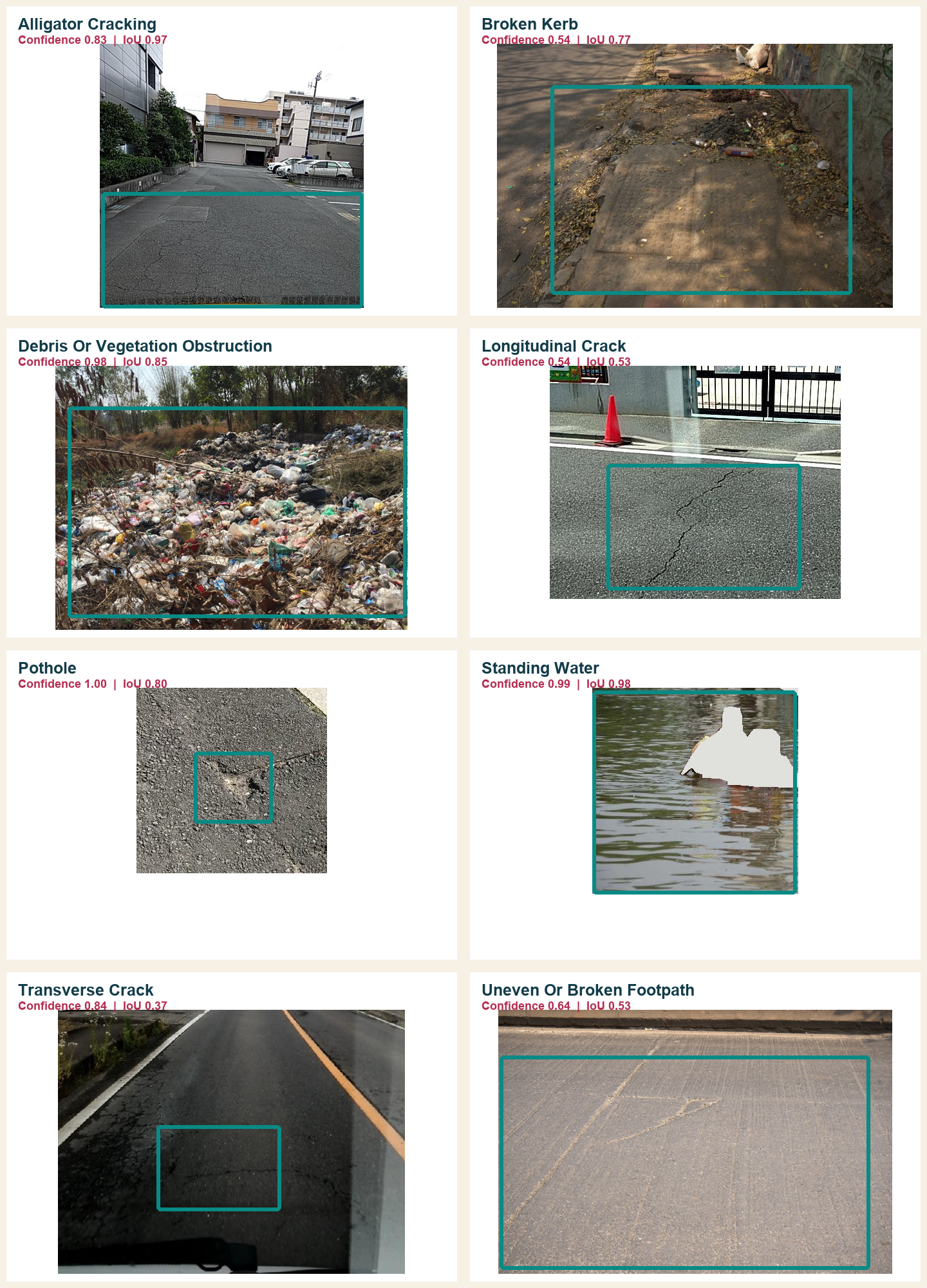}
  \caption{Successful held-out detections across the eight-class condition taxonomy. Each panel reports model confidence and intersection over union for the localised condition.}
  \label{fig:success}
\end{figure}

Figure~\ref{fig:success} foregrounds the strongest correct held-out examples for all eight visible-condition classes. The gallery includes road cracks, potholes, standing water, footpath damage, broken kerbs and obstruction. Several examples combine high confidence with close spatial agreement. Standing water reaches 0.99 confidence and 0.98 intersection over union. Alligator cracking reaches 0.83 confidence and 0.97 intersection over union. The examples demonstrate that one shared perception framework can represent road, walking and kerb conditions instead of restricting inspection to pavement defects.

\FloatBarrier

\begin{figure}[H]
  \centering
  \begin{subfigure}[t]{0.96\textwidth}
    \centering
    \includegraphics[width=0.86\linewidth]{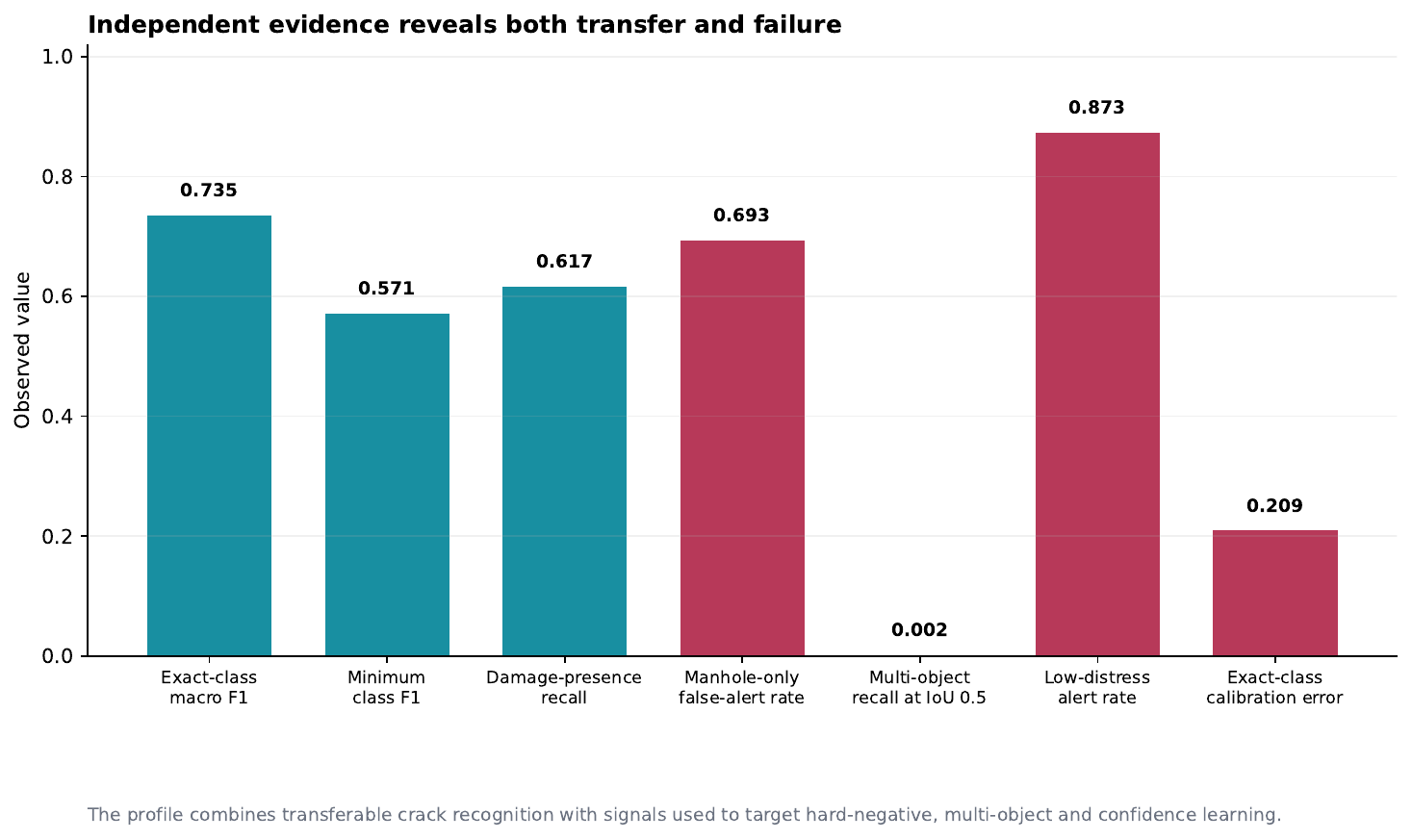}
    \caption{Cross-source learning profile. Strong crack-class transfer is considered with the signals used to target the next learning cycle.}
    \label{fig:validation}
  \end{subfigure}

  \vspace{4pt}
  \begin{subfigure}[t]{0.65\textwidth}
    \centering
    \includegraphics[width=\linewidth]{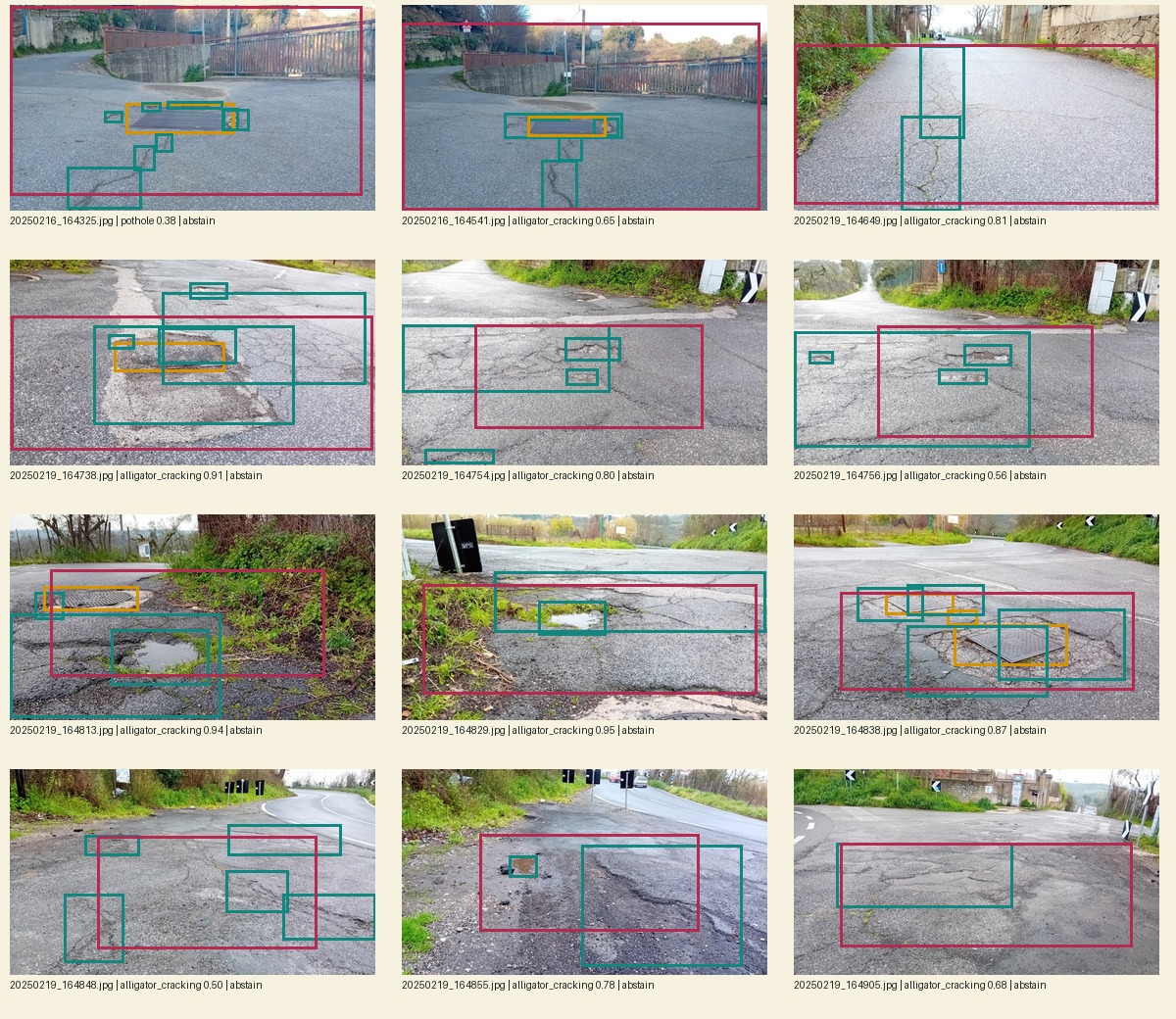}
    \caption{Compact hard-case review. Teal marks reviewed damage, amber marks manhole distractors and red shows model predictions.}
  \end{subfigure}
  \caption{Cross-source recognition evidence. The transfer profile remains the primary result, while a smaller diagnostic panel shows how difficult scenes guide hard-negative and multi-object learning.}
  \label{fig:detectionevidence}
\end{figure}

Figure~\ref{fig:detectionevidence} makes these abstract measures visible through held-out and cross-source examples.

Model cards and dataset records report the selected model, intended use, limits and provenance \citep{mitchell2019model,gebru2021datasheets}. The next learning cycle expands multi-object localisation, adds hard negatives and deepens local transfer. This creates a clear collaboration pathway for agencies and researchers that can contribute lawful regional imagery, expert review or deployment settings.

\subsection{From perception to decision support}

The purpose of condition detection is not to produce boxes for their own sake. The longer-term workflow asks questions such as ``Which footpaths should be repaired first?'' A useful answer should consider defect type, severity, confidence, human review, network continuity, access to schools, hospitals and public transport, mobility needs, flood history, cost and the expected benefit of repair.

\begin{figure}[H]
  \centering
  \includegraphics[width=0.94\textwidth]{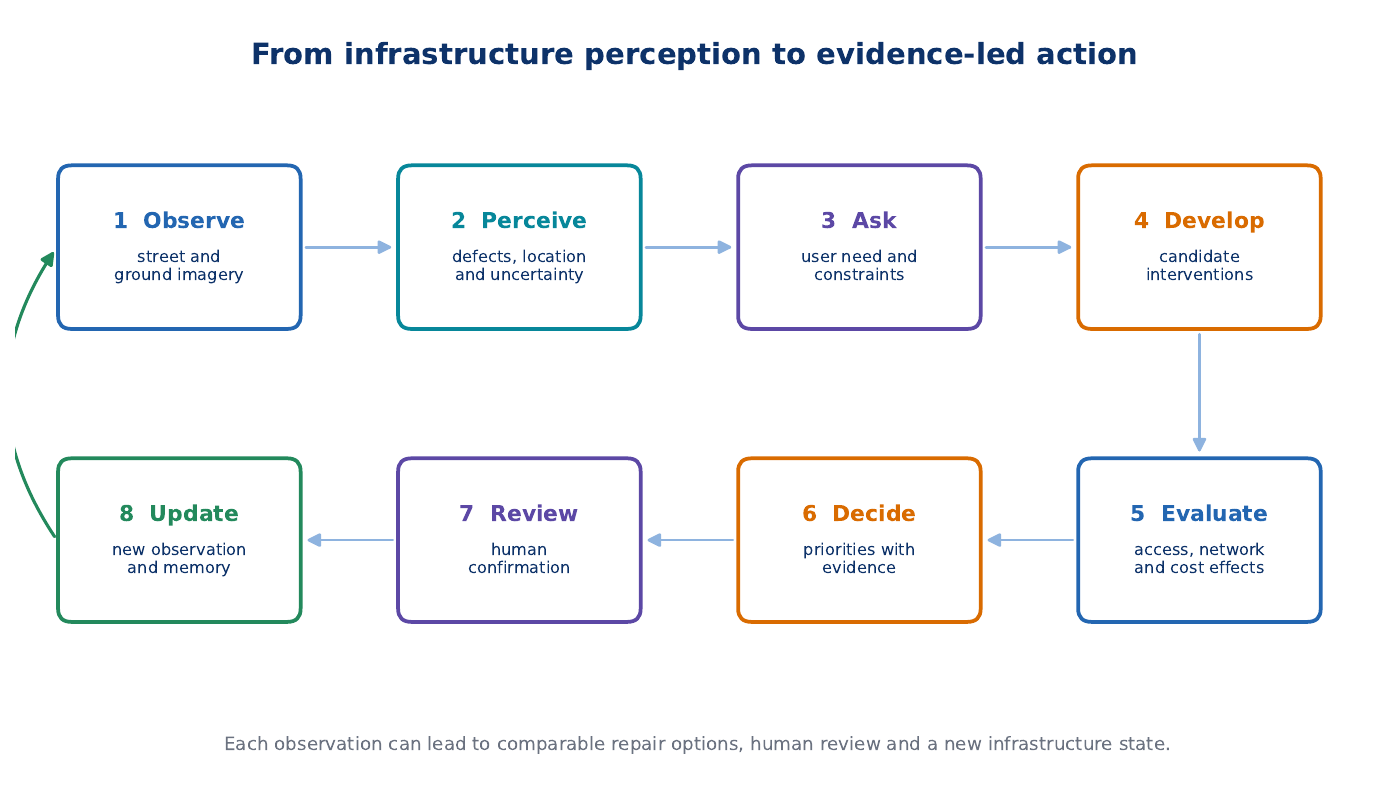}
  \caption{Perception-to-decision feedback loop. Implemented perception and analysis functions form the foundation. Automated intervention optimisation and closed-loop field validation are collaboration priorities.}
  \label{fig:loop}
\end{figure}

Figure~\ref{fig:loop} shows the proposed closed loop. Observation and perception create a current infrastructure state. A user question defines the decision. Candidate interventions are developed and evaluated through accessibility, network and cost models. A human reviews the recommendation. New imagery then checks whether the intervention changed the observed condition.

This workflow resembles a digital twin because it connects observations with a changing digital state. It becomes closer to a world model when the system can compare candidate actions, estimate consequences, retain evidence and learn from feedback. Footpath decisions need not use MATSim. Many can use pedestrian-network connectivity, accessible-route analysis, service-population measures and budget optimisation. Agent-based walking simulation becomes useful when crowd behaviour or detailed interactions matter.

\section{One platform, several ways to work}
\label{sec:using-platform}

\begin{figure}[H]
  \centering
  \includegraphics[width=0.98\textwidth]{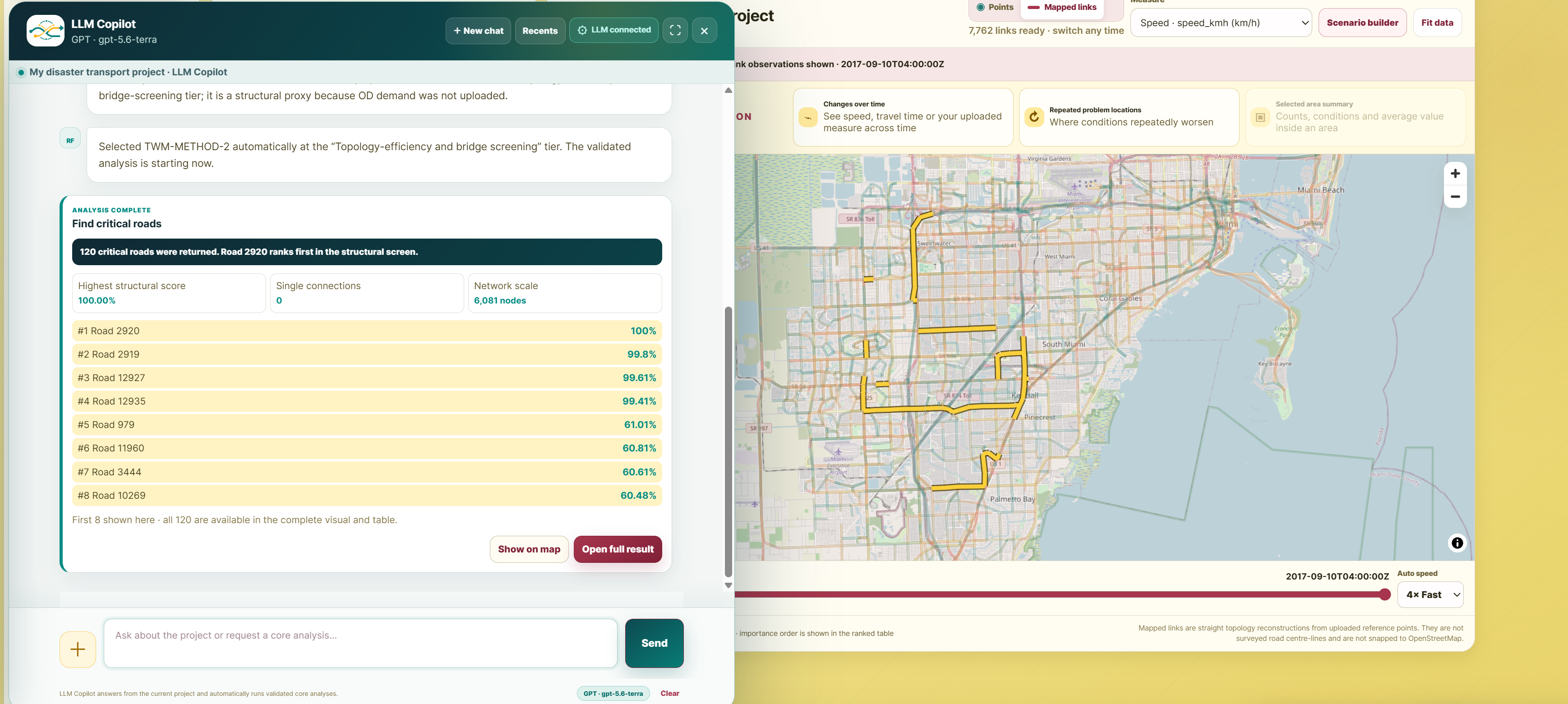}
  \caption{Question-led critical-road analysis. The LLM Copilot selects the documented analysis tier, presents the ranked result and connects it to the roads displayed on the geographical workspace.}
  \label{fig:llmcriticalroads}
\end{figure}

ResiliFlow is designed for people who begin with different evidence and different questions. A transport analyst can work directly with mapped entities and equations. An infrastructure manager can begin with a street view. An emergency planner can compare disruption scenarios. A researcher can inspect the full validation trail. The local Assistant supports navigation, while the optional LLM Copilot can prepare and execute bounded calls to the six documented tools. Figure~\ref{fig:llmcriticalroads} shows the Copilot selecting a critical-road method, returning a ranked result and displaying the corresponding roads on the shared map. The executable network module performs the calculation, while the Copilot manages the question, inputs and explanation.

\begin{figure}[H]
  \centering
  \includegraphics[width=0.94\textwidth]{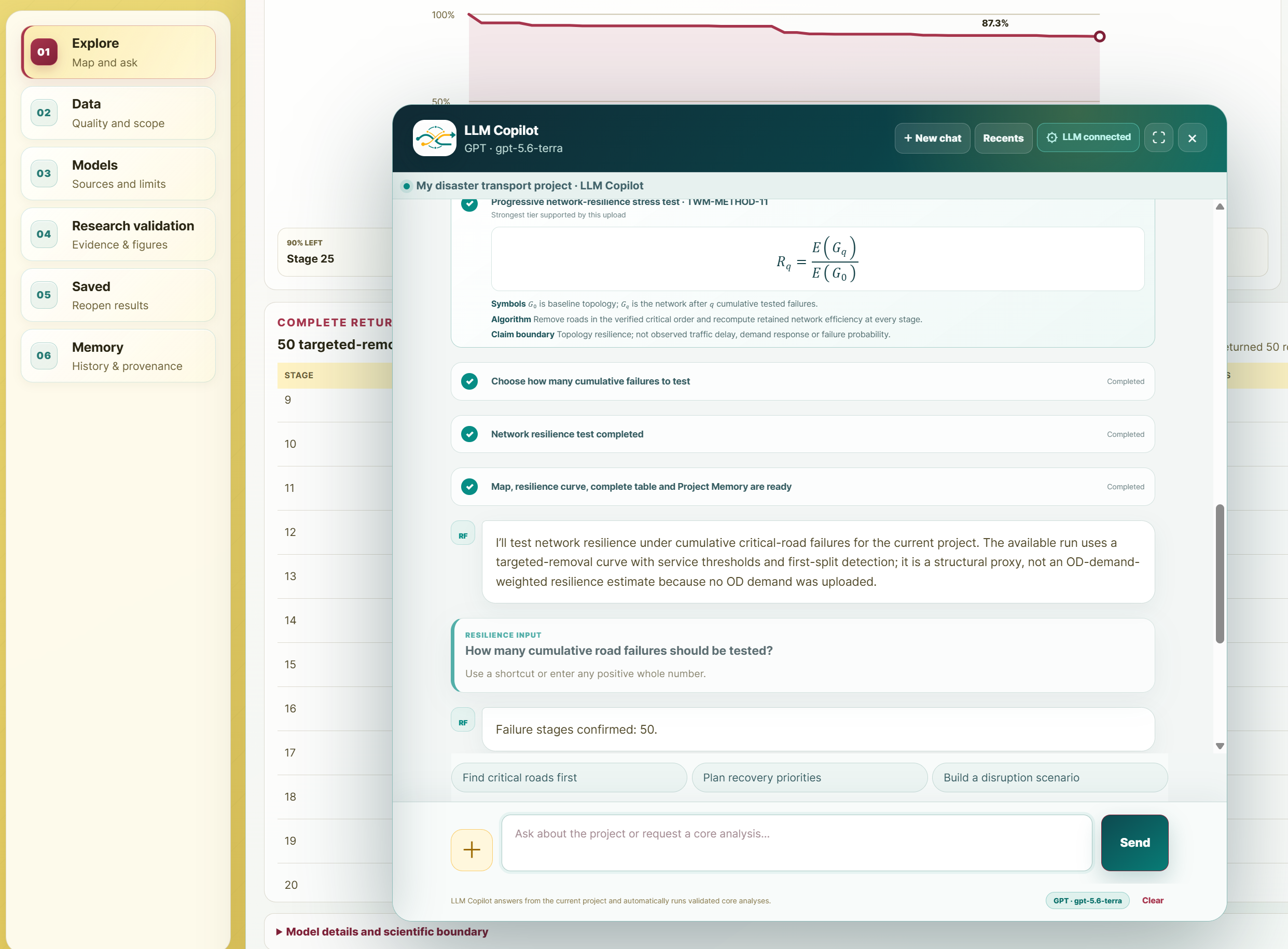}
  \caption{LLM-guided network-resilience testing. The interface exposes the model definition, gathers the requested number of cumulative failures and links the completed run to its curve, table and project record.}
  \label{fig:llmresilience}
\end{figure}

The same interaction pattern supports multi-step analyses. Figure~\ref{fig:llmresilience} shows a network-resilience test in which the Copilot exposes the performance definition, confirms the cumulative-failure input, runs the executable test and records completion. The result remains connected to the resilience curve, return table and Project Memory behind the conversation.

\begin{figure}[H]
  \centering
  \includegraphics[width=0.98\textwidth]{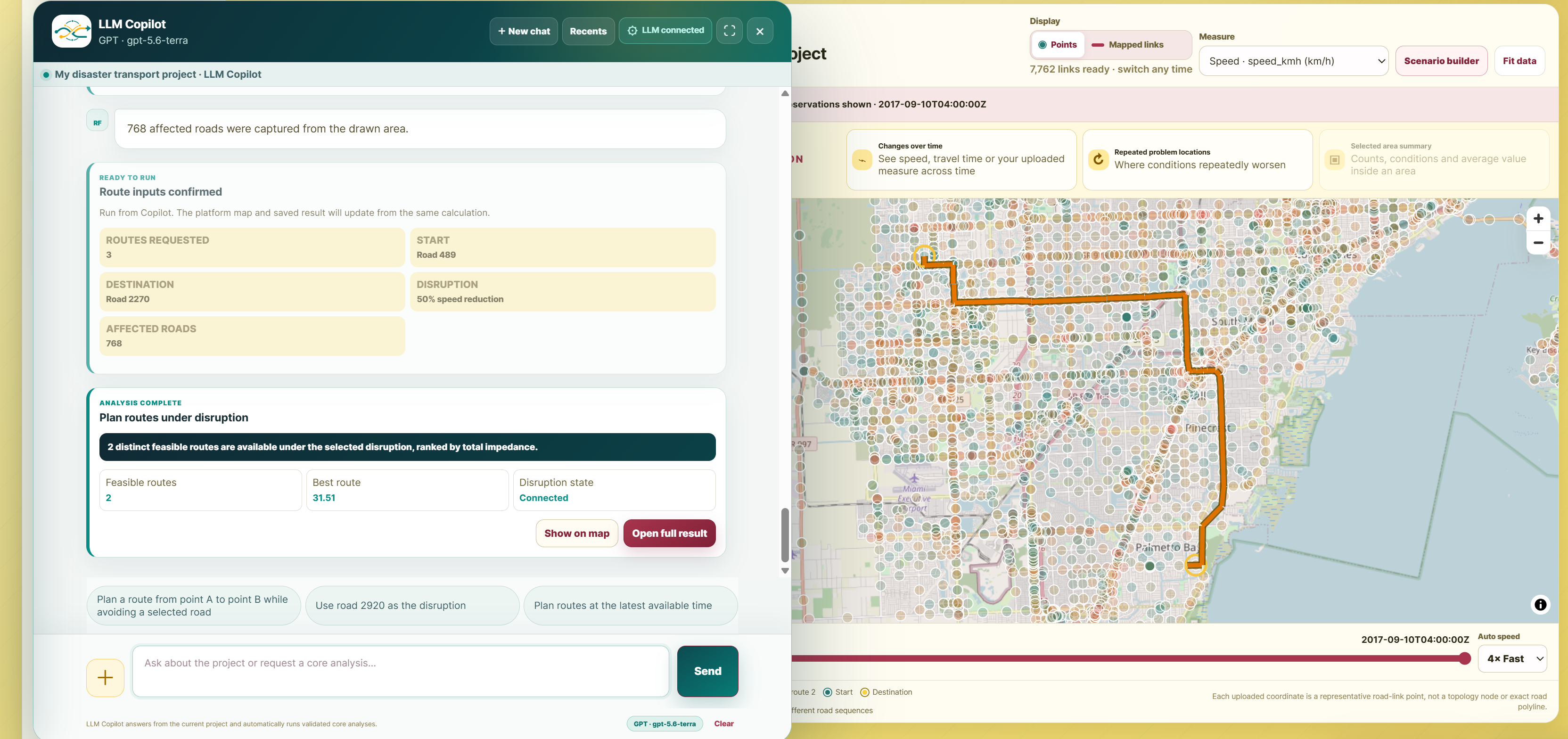}
  \caption{LLM-guided routing under disruption. Confirmed route inputs and affected roads are passed to the routing module, which returns feasible alternatives and displays the selected path on the shared map.}
  \label{fig:llmrouting}
\end{figure}

Question-led operation also supports spatial selections and routing under a specified disruption. In Figure~\ref{fig:llmrouting}, the Copilot confirms the affected roads, origin, destination, route count and speed reduction before execution. It then returns feasible alternatives and displays the selected route on the same map. This creates a visible path from a user's request to model inputs, computation and an actionable geographical result.

\FloatBarrier

The two workspaces form one extensible transport world model. Disaster analysis explains network consequences and response options. Infrastructure perception adds observed conditions. Research Validation makes model reasoning visible. The Assistant and optional Copilot help users reach these capabilities. Memory retains the project state needed to revisit or extend an analysis.

\section{Technical foundations and extensibility}

\subsection{A modular model library}

The platform is designed around model contracts. A contract states required inputs, assumptions, parameters, outputs, validation rules and failure states. It allows a graph metric, optimiser, detector or simulation engine to participate in the same cycle. This approach supports method comparison while preserving a stable user workflow.

The existing research portfolio demonstrates the range of models that such a library can connect. Risk and cascading-failure models capture multimodal dependence \citep{guo2021risk,guo2020risk}. Evacuation models address behaviour, shelter assignment and traffic dynamics \citep{ng2010a,ng2010reliable,dixit2012modeling,dixit2014evacuation}. Network design models include behavioural choice, uncertain demand and reliability \citep{xu2024exploring,xu2024predisaster,xu2025predisaster}. Economic models connect transport and wider production systems \citep{robson2017a,robson2018a,shahriari2023integrating,wang2024calibration}. AI models add visible infrastructure state \citep{arya2024rdd,li2024rddyolo,zhao2024detrs,jocher2026yolo26}.

Transport investment research also contributes decision logic. \citet{li2015transit} studied technology selection under population volatility. \citet{guo2018how} examined time-inconsistent preferences in rail investment. \citet{guo2018stochastic} developed real-option switching between fixed and flexible transit. \citet{sun2017evolution} analysed the evolution of transit modes. \citet{sun2018optimal} studied rail-line extension, and \citet{guo2021timedependent} optimised time-dependent fares. \citet{guo2023investment} considered investment timing and line length under demand uncertainty, while \citet{li2025joint} jointly optimised airport-terminal expansion timing and increment. These methods show how perception and disruption evidence can later feed investment and recovery decisions.

\subsection{Model choice as evidence}

A world model should not hide model selection. If several algorithms can answer a question, the platform records why one was chosen. The choice may depend on data, scale, runtime, interpretability and validation. For example, a shortest-path method is appropriate for a single routing query, while an equilibrium assignment is needed when many travellers affect each other's costs. A compact detector may suit edge deployment, while a larger transformer may improve spatial accuracy. The decision record turns these trade-offs into evidence.

The model library also makes disagreement useful. Two models can be run on the same state, and their output difference can be reported. Cross-model comparison is especially valuable when behaviour and demand are uncertain. It can reveal which conclusions are stable and which depend on an assumption. This is more informative than presenting one preferred output without context.

\subsection{Safe openness}

``Open'' describes the intended collaboration model. It does not mean that every component or third-party data source can be redistributed without review. The platform is proceeding through UNSW intellectual property and legal due diligence under technology ID 2026-179. The planned release will separate open platform code and examples from user secrets, provider keys, restricted data, model weights that need licence review and institution-specific deployment settings.

The release structure is intended to let contributors add a model through a documented contract, reproduce a published example and inspect an evidence record. Sensitive operational details are omitted from this paper. This balance supports scrutiny and contribution while protecting user data and third-party terms.

\section{Future opportunities}

ResiliFlow already brings executable network analysis, model visualisation, user guidance, project memory and infrastructure perception into one platform. Its open architecture creates several immediate collaboration opportunities. Regional partners can contribute lawful road and footpath imagery. Transport authorities can define priority questions and test repair workflows. Modelling groups can add transparent simulation, accessibility, economics and optimisation modules. Human-centred researchers can evaluate how practitioners understand confidence, trade-offs and model explanations.

The next growth cycle focuses on four connected capabilities. The first is broader perception across varied road environments. The second links confirmed conditions to network continuity, access, equity, risk and repair cost. The third closes the feedback loop through new observations after an intervention. The fourth expands the model library through reproducible contracts and shared examples. Satellite hazard context will add a regional view of flooding, land movement and disruption, while street-level perception will retain the detail needed for road, footpath and kerb management.

Broader research can connect shared autonomous vehicles, transit investment, equity, evacuation and economic impact. \citet{sevim2025a} showed how shared autonomous vehicles can support rural evacuation. \citet{li2026vehicle} used deep learning for automated-guideway-transit occupancy. \citet{waller2025mobility} framed mobility as a resource for human-centred automation. Safety research on automated driving and risk attitudes adds further behavioural evidence \citep{dixit2019risk,arbis2016impact}. These directions can be added as interoperable modules when their evidence and user purpose are clear.

The most important future test is whether the platform improves real work. That requires partnerships with transport authorities, infrastructure managers, emergency planners and community organisations. Useful evaluation will measure task completion, decision traceability, uncertainty comprehension, time saved and whether recommended actions remain fair across locations and population groups.

\section{Conclusions}

ResiliFlow presents a transport world model as an executable and inspectable cycle. It connects observations, established transport models, learned perception, validation, decisions, feedback and reusable memory. This definition gives the world model idea a practical meaning for infrastructure resilience, response and recovery.

The implemented platform already joins two substantial workspaces. Disaster Transport Resilience Analysis provides six connected tools, Research Validation, optional simulation, local guidance, an LLM Copilot and project memory. AI-based Transport Infrastructure Perception and Decision Support brings street-level and satellite context into a governed evidence workflow for road, footpath and kerb conditions.

The perception results show that a common detector can recognise eight visible-condition classes across road, footpath and kerb evidence. Cross-source learning, hard-case analysis and the Bhopal interface create a direct path from image evidence to mapped conditions and repair options.

The next step is open, governed collaboration. ResiliFlow offers a common structure in which new datasets, models, simulations and decision methods can be contributed without losing provenance or human responsibility. Its central promise is simple: transport data should lead to understandable evidence, and understandable evidence should lead to better tested action.

\section*{Author contributions}

\noindent
\textbf{Junxiang Xu:}
Conceived the Transport World Model platform, designed its system architecture and led the complete prototype development. Developed the platform software, transport analysis functions, simulation workflows, model library, LLM orchestration, Research Validation workflow and infrastructure perception components. Conducted data processing, functional testing, computational experiments, system validation, interface development and results visualisation. Led project administration and prepared the working paper.

\noindent
\textbf{Vinayak Dixit:}
Co-conceived the platform and contributed to its methodological design, technical development strategy, functional testing and system validation. Supervised the project and provided continuing technical guidance. Reviewed the working paper.

\noindent
\textbf{S. Travis Waller:}
Co-conceived the platform and contributed to its methodological foundations, technical direction and system validation. Provided strategic technical guidance and reviewed the working paper.

\noindent
\textbf{Divya Jayakumar Nair:}
Provided project resources, contributed to platform testing and functional validation, and supported external collaboration, partner engagement and coordination. Reviewed the working paper.

\noindent
\textbf{Qianwen (Vivian) Guo:}
Provided data resources, supported data preparation and contributed to platform testing and validation. Reviewed the working paper

\noindent
\textbf{Sisi Jian:}
Provided data resources, supported data preparation and contributed to platform testing and validation. Reviewed the working paper

\noindent
\textbf{Xiao Wen:}
Provided data resources, supported data preparation and contributed to platform testing and validation. Reviewed the working paper

\noindent
\textbf{Ashutosh Ashutosh:}
Provided technical advice, platform feedback and functional validation. Reviewed the working paper

\noindent
\textbf{Sunhyung Yoo:}
Provided technical advice, platform feedback and functional validation. Reviewed the working paper

\noindent
\textbf{Julius Secadiningrat:}
Provided technical advice, platform feedback and functional validation. Reviewed the working paper

\noindent
\textbf{Jingni Guo:}
Contributed to user-interface implementation, visual design, platform styling and page layout. Reviewed the working paper

\section*{Code availability}

Selected design code and components of the prototype testing platform will be progressively released at \url{https://github.com/ResiliFlow}, subject to UNSW intellectual property and legal review. 

\begingroup
\onehalfspacing
\setlength{\bibsep}{0.5\baselineskip}
\sloppy
\bibliographystyle{elsarticle-harv}
\bibliography{references}
\endgroup
\end{document}